\documentclass[aap]{imsart}

\RequirePackage{amsthm,amsmath,amsfonts,amssymb}
\RequirePackage[numbers]{natbib}
\RequirePackage[colorlinks,citecolor=blue,urlcolor=blue]{hyperref}
\usepackage{latexsym}
\usepackage{amsbsy}
\usepackage{pstricks}
\usepackage{pst-fill,pst-grad}
\usepackage{color}
\allowdisplaybreaks

\startlocaldefs

\numberwithin{equation}{section}
\newtheorem{theorem}{Theorem}[section]
\newtheorem{corollary}[theorem]{Corollary}
\newtheorem{lemma}[theorem]{Lemma}
\newtheorem{proposition}[theorem]{Proposition}

\theoremstyle{definition}

\newtheorem{assumption}[theorem]{Assumption}
\newtheorem{remark}[theorem]{Remark}

\newtheorem{example}[theorem]{Example}

\newdimen\AAdi%
\newbox\AAbo%
\def\AAk#1#2{\setbox\AAbo=\hbox{#2}\AAdi=\wd\AAbo\kern#1\AAdi{}}%

\makeatletter
\def\eqlabel#1{\def\@currentlabel{#1}}

\def\formula#1{\def\@tempa{#1}\let\@tempb\theequation\def\theequation{%
\hbox{#1}}\def\@currentlabel{(\theequation)}$$}
\def\endformula{\leqno\hbox{(\@tempa)}$$\@ignoretrue\let\theequation\@tempb}

\def\given{\hskip5\p@\relax\vrule\@width.4\p@\hskip5\p@\relax}

\newcommand{\open}[1]{%
\par\normalfont\topsep6\p@\@plus6\p@\trivlist\item[\hskip\labelsep\itshape#1%
\@addpunct{.}]\ignorespaces}

\DeclareRobustCommand{\close}[1]{%
  \ifmmode 
  \else \leavevmode\unskip\penalty9999 \hbox{}\nobreak\hfill
  \fi
  \quad\hbox{$#1$}}
\makeatother

\newlength{\toskip}

\def\<{\langle}
\def\>{\rangle}

\def \Var {\textrm{Var}}
\def \Ent {\textrm{Ent}}
\def \Osc {\textrm{Osc}}
\def \Cov {\textrm{Cov}}

\endlocaldefs

\begin{document}

\begin{frontmatter}

\title{Diffusion annealed Langevin dynamics: a theoretical study.}

\runtitle{Annealed Langevin}

\begin{aug}

\author[A]{\fnms{Patrick}~\snm{Cattiaux}\ead[label=e1]{patrick.cattiaux@math.univ-toulouse.fr}},
\author[B]{\fnms{Paula}~\snm{Cordero-Encinar}\ead[label=e2]{paula.cordero-encinar22@imperial.ac.uk}}
\and
\author[C]{\fnms{Arnaud}~\snm{Guillin}\ead[label=e3]{arnaud.guillin@uca.fr}}

\address[A]{Institut de Math\'ematiques de Toulouse. CNRS UMR 5219. \\
Universit\'e de Toulouse, 118 route
de Narbonne, F-31062 Toulouse cedex 09.\\ \printead[presep={\ }]{e1}}

\address[B]{Department of Mathematics. Imperial College London. London, UK.\\
\printead[presep={\ }]{e2}}

\address[C]{Universit\'e Clermont Auvergne, CNRS UMR 6620, LMBP, F-63000 Clermont-Ferrand, France.\\
\printead[presep={\ }]{e3}}
\end{aug}

\begin{abstract}
In this work we study the diffusion annealed Langevin dynamics, a score-based diffusion process recently introduced in the theory of generative models and which is an alternative to the classical overdamped Langevin diffusion. Our goal is to provide a rigorous construction and to study the theoretical efficiency of these models for general base distribution as well as target distribution. As a matter of fact these diffusion processes are a particular case of Nelson processes i.e. diffusion processes with a given flow of time marginals.

Providing existence and uniqueness of the solution to the annealed Langevin diffusion leads to proving a Poincar\'e inequality for the conditional distribution of $X$ knowing $X+Z=y$ uniformly in $y$, as recently observed by one of us and her coauthors.
Part of this work is thus devoted to the study of such Poincar\'e inequalities. Additionally we show that strengthening the Poincar\'e inequality into a logarithmic Sobolev inequality improves the efficiency of the model.
\end{abstract}



\end{frontmatter}

\section{Introduction}\label{secintro}

The aim of this paper is to give a rigorous presentation of the recently introduced \emph{diffusion annealed Langevin dynamics} \cite{conf6}. This stochastic process is a score based generative model and provides an alternative to the well known overdamped Langevin process  and its reversed in time version commonly used for sampling purpose. In particular, we will fill some gaps in the main arguments used for building the annealed Langevin dynamics discussed in \cite{conf6,Guo,Corderoetal}. We will not discuss its practical efficiency nor its numerical counterparts, that is we will not introduce nor discuss the corresponding discrete algorithms, presented in \cite{Corderoetal} by the second author, and the references therein. However, some quantitative aspects, useful for discretization schemes or important from the statistical point of view, are discussed in details. 
 
 Also, for distributions like the gaussian, an important idea introduced in the papers on diffusion annealed Langevin dynamics  consists in using a functional inequality (namely the Poincar\'e inequality) to control some covariance. This inequality is crucial in \cite{Corderoetal} for proving that the score of the intermediate distributions is Lipschitz continuous, which, as we explain in Section~\ref{secstoc}, ensures the existence and uniqueness of strong solutions for the annealed Langevin diffusion. As a matter of fact, heavy tailed base distributions are also particularly well suited for the model as will see in an example.

More generally, Section \ref{secstoc} provides a rigorous introduction to annealed Langevin dynamics, and concludes in the main theorem, which establishes conditions for the strong and weak existence and uniqueness of solutions, as well as a quantitative estimate of the Kullback–Leibler divergence between the annealed Langevin dynamics and the flow of time marginals targeted by the annealed Langevin diffusion, i.e. a neat control of the bias of the algorithm.

We then show that a Poincaré inequality for the conditional distribution is necessary to guarantee existence and uniqueness of solutions.
In doing so, we complete the results in \cite{Corderoetal}, by recalling and applying several known results, concerning the Poincar\'e constant of perturbed measures, as explained in \cite{CGperturb}, see Sections \ref{secgame} and \ref{secaboutpoinc}. We then prove new results for this Poincar\'e constant, Section \ref{secnewpoinc}. 

Finally, in Section \ref{seclogsob}, we demonstrate that one can improve the bias between the target distribution and the final distribution of the annealed Langevin dynamics by using the stronger logarithmic Sobolev inequality.
\medskip

\section{Score based Stochastic dynamics}\label{secstoc}

\subsection{Generative diffusion model}
A central point in modern Statistics is to learn a probability distribution $\pi$ from a given set of samples. Recent advances proposed to introduce dynamical systems bridging the target distribution $\pi$ and a base (or noising) distribution $\nu$. Several methods have been used, see e.g. \cite{AVE}, in particular using stochastic differential equations.

We only consider the case where $\mathbb R^d$ is the state space. A popular method is to build an ergodic diffusion process with invariant distribution given by $\nu$ and initial distribution $\pi$ (overdamped Langevin which can be implemented using the Langevin-Monte Carlo algorithm). If $\nu(dx)= e^{-V(x)} dx$ for some smooth enough $V$ (say $C^2$), the natural associated process is the solution of 
\begin{equation}\label{eqLMC}
dX_t=\sqrt 2 \, dB_t \, - \, \nabla V(X_t) \, dt \quad \textrm{ with } \quad \mathcal L(X_0)=\pi \, .
\end{equation}
where $B_.$ is a standard Brownian motion. Since $\nabla V$ is local Lipschitz, existence and uniqueness hold true up to an explosion time. In order to guarantee that the explosion time is infinite some additional assumptions are required like for instance
\begin{enumerate}
\item [(H1)] \quad there exists some $\psi$ such that $\psi(x) \to +\infty$ as $|x| \to +\infty$ and $\Delta \psi - \nabla V.\nabla \psi$ is bounded above,
\item [(H2)] \quad $\int |\nabla V|^2 d\nu < +\infty$. 
\end{enumerate}
(H1) is immediate using It\^{o} formula, (H2) is more intricate and requires the use of Dirichlet forms and some tools in Potential Theory (see e.g. \cite{CF96}). In these cases $\nu$ is a reversible (hence invariant) measure for the dynamics \eqref{eqLMC}.

Using the time reversed process, one is able to rebuild the target measure $\pi$. Unfortunately, the invariant measure is attained in an infinite time (a possible method to overcome this difficulty is to consider bridges, see \cite{Leo-survey} for a survey). Hence one is obliged to stop the procedure at a given time $T$, to evaluate the bias introduced by this cut-off using an appropriate distance between measures, and then to study the time reversed diffusion. 

More precisely, if $p_t(.)$ denotes the density of $X_t$ w.r.t. the Lebesgue measure (which exists thanks to ellipticity and regularity), the time reversed process $X^T_t=X_{T-t}$ would satisfy
\begin{equation}\label{eqLMCrev}
dX^T_t=\sqrt 2 \, dB_t \, + \, \nabla V(X^T_t) \, dt \, + \, 2 \, \nabla.\ln p_{T-t}(X^T_t) \, dt
\end{equation}
with initial distribution $\mathcal L(X^T_0)=p_T(x) dx$. Such time reversal argument requires to be careful with the assumptions. One possible set of assumptions is that $\nabla V$ is (global) Lipschitz and $\nabla.p_.(.) \in \mathbb L^1_{loc}(dt,dx)$. For this result see e.g. \cite{MNS} Theorem 2.3. Other sets of assumptions are possible. We will not discuss them here but in another paper.

Unfortunately, there is no explicit method to recover the density $p_.$ (or the gradient of its logarithm $\nabla. \ln p_.$ called the score in the statistical literature), except in the particular case where $\nu$ is a gaussian measure. Assume for simplicity that $\nu$ is the standard normal distribution, then $X_.$ is the classical Ornstein-Uhlenbeck process, starting from $\pi$, so that 
\begin{equation}\label{eqOUconvol}
X_t = e^{-t} X_0 + \sqrt{1-e^{-2t}} \, G
\end{equation}
where $G$ is a standard normal distribution independent of $X_0$, implying that the law at time $t$ is given by
\begin{equation}\label{eqOUdens}
p_t(x) \, dx = \frac{\pi(x/e^{-t})}{e^{-td/2}}*\frac{\nu(x/(1-e^{2t})^{\frac 12})}{(1-e^{-2t})^{d/2}} \, .
\end{equation}
The remaining problem is that at the stopping time $T$, $X_T$ is not equal to the gaussian distribution. However it is not far from it. Indeed, if we assume that the Kullback-Leibler divergence (or relative entropy) $$d_{KL}(\pi,\nu): = \int \, \ln(d\pi/d\nu) \, d\pi$$ is finite, then 
\begin{equation}\label{eqentropgauss}
d_{KL}(\pi_T,\nu) \, \leq \, e^{-2T} \, d_{KL}(\pi,\nu) \, ,
\end{equation} 
since the gaussian distribution $\nu$ satisfies a logarithmic Sobolev inequality $$\Ent_\nu(f) := \int \, f \, \ln\left(\frac{f}{\int f d\nu}\right) \, d\nu \, \leq \frac 12 \, \int \frac{|\nabla f|^2}{f} \, d\nu \, ,$$ for any smooth and positive $f$.

If the SDE  
\begin{equation}\label{eqLMCrev2}
dZ_t=\sqrt 2 \, dB_t \, + \, Z_t \, dt \, + \, 2 \, \nabla.\ln p_{T-t}(Z_t) \, dt
\end{equation}
 with initial distribution $\mathcal L(Z_0)=\nu$ admits a unique strong solution, one has 
 \begin{align*}
     d_{KL}(\pi,\mathcal L(Z_T)) &= d_{KL}(\mathcal L(X_T^T),\mathcal L(Z_T)) \leq d_{KL}(\mathcal L(X_0^T),\mathcal L(Z_0)) \\
     &= d_{KL}(\pi_T,\nu) \leq e^{-2T} \, d_{KL}(\pi,\nu) \, .
 \end{align*} In order to ensure that \eqref{eqLMCrev2} has a unique solution, the most standard way is to get conditions for $\nabla.\ln p_{T-t}$ to be Lipschitz in $x$ uniformly in $t\in[0,T-\varepsilon]$ for any $\varepsilon >0$ (we shall say uniformly on $[0,T^-)$), in other words conditions for $\nabla^2 \ln p_t$ to be bounded, uniformly in time on $(0^+,T]$. For a unique weak solution it is enough to have that $\nabla.\ln p_{T-t}$ is bounded in $(x,t)$ for $t\in[0,T-\varepsilon]$ for any $\varepsilon >0$.
\medskip

\subsection{Annealed Langevin process}
A natural alternative way is to introduce some interpolation between $\pi$ and $\nu$. Several interpolation procedures are possible, see e.g. \cite{Albergo1}. From the probabilistic point of view a very natural one is the following analogue of what we have done before in the gaussian case: let $X$ and $Z$ be two independent random variables with respective distributions $\pi$ and $\nu$, consider the interpolation flow 
\begin{equation}\label{eqinterpolationgeneral}
\tilde X_t= \sqrt{\lambda_t} \, X + \sqrt{1-\lambda_t} \, Z
\end{equation} 
for an increasing  flow $t\in[0,T] \mapsto \lambda_t$ of parameters in $[0,1]$ (called the schedule in the machine learning literature) such that hopefully $\lambda_T=1$ and  $\lambda_0=0$. We denote $\tilde \pi$ (resp. $\tilde \nu$) the distribution of $\tilde X_T$ (resp. $\tilde X_0$).

At time $t$ the process has a distribution given by 
\begin{equation*}\label{eqconvoldens}
p_t(x) \, dx = \frac{\pi(x/\sqrt{\lambda_t})}{\lambda_t^{d/2}}*\frac{\nu(x/\sqrt{1-\lambda_t})}{(1-\lambda_t)^{d/2}} \, .
\end{equation*} 
The question is then to represent $X_.$ as the solution of some Stochastic Differential Equation.
\medskip

For this purpose, the first thing to do is to show that $p_.$ does satisfy a Fokker Planck equation. For convenience, one introduces a renormalization parameter $\kappa$ and the time-changed flow 
\begin{equation}\label{eqhat}
\hat p_t=p_{\kappa t} \quad \textrm{ for $t \in [0,T/\kappa]$} \, . 
\end{equation}
Inspired by optimal transport, one way consists in finding a vector field $v_.$ satisfying the transport (continuity) equation 
\begin{equation}\label{eqtransport}
\partial_t \hat p_t + \nabla.(v_t \, \hat p_t) =0
\end{equation}
 so that, at least formally,
\begin{equation}\label{eqFPannealed}
\partial_t \hat p_t = \Delta \hat p_t \, - \, \nabla.((\nabla \ln \hat p_t +v_t)\hat p_t) \, .
\end{equation}
To this end one can for instance apply Theorem 8.3.1 in \cite{AGS}, i.e
\begin{theorem}\label{thmambrosio}
Assume that for all $t \in [0,T]$, $\int |x|^2 \, p_t(x) \, dx < +\infty$. Define, when it exists, $$|p'_t| = \lim_{h \to 0} \; \frac{1}{|h|} \; W_2(p_{t+h},p_t)$$ where $W_2$ denotes the $2$-Wasserstein distance. Assume that $|p'_.| \in \mathbb L^1([0,T])$. Remark that $|(\hat p)'_t|=\kappa \, |p'_{\kappa t}|$ for all $t \in [0,T/\kappa]$.

Then there exists a Borel vector field $v_.$ such that \eqref{eqtransport} is satisfied and $$||v_t||_{\mathbb L^2(\hat p_t dx)} \, = \, \kappa \, |p'_{\kappa t}| \quad \textrm{ for almost all } \; t \in [0,T/\kappa] \, .$$
\end{theorem}

This result is recalled and used in \cite{Guo,Corderoetal}. The next step in these papers is to introduce first the so called \emph{diffusion annealed Langevin dynamics} given by the S.D.E. on $[0,T/\kappa]$,
\begin{equation}\label{eqannealed}
dY_t = \sqrt 2 \, dB_t + \nabla \ln \hat p_{t}(Y_t) \, dt \quad \textrm{ with } \quad \mathcal L(Y_0)= \tilde \nu \, ,
\end{equation}
which, in a sense, is mimiking \eqref{eqLMCrev}. If the drift is explicit, the law at time $t$ is not, once again. To ensure that the law at time $t$ is given by $\hat p_t$ one has to consider
\begin{equation}\label{eqannealedbis}
dX_t = \sqrt 2 \, dB_t + (\nabla \ln \hat p_{ t} + v_t)(X_t) \, dt \quad \textrm{ with } \quad \mathcal L(X_0)= \tilde \nu \, .
\end{equation}
This time the drift is not explicit but the law at time $t$ should be $\hat p_t$. 

The first problem is to show existence and, if necessary uniqueness, of the solution of \eqref{eqannealed} and \eqref{eqannealedbis}. One also has to check that the marginals flow of a (or the) solution of \eqref{eqannealedbis} is actually $\hat p_.$. The next problem is to use Girsanov theory in order to control the Kullback-Leibler divergence between the distributions on the path space. To this end, uniqueness for \eqref{eqannealed} is required. These problems are not completely discussed in \cite{Guo,Corderoetal} and some arguments are missing. 

We shall thus fill the gap(s) and complete the proofs of the main results in these references.
\medskip

Existence and strong uniqueness for \eqref{eqannealed} are ensured as soon as $\nabla^2 \ln \hat p_{ .}$ is bounded uniformly in $t$ for $t \in [0,(T/\kappa)^-]$. If $\nabla \ln \hat p_{ .}$ is bounded, one may use Girsanov theorem and get existence and weak uniqueness.

\eqref{eqannealedbis} is much more delicate. Indeed the only known property on $v_t$ is its square integrability w.r.t. $\hat p_t dx$, so that none of the standard existence theorem applies. If in addition to the uniqueness of the solution of \eqref{eqannealed}, we assume that $$\int_0^{(T/\kappa)} \, ||v_t||^2_{\mathbb L^2(\hat p_t dx)} dt < + \infty \, ,$$ i.e. the finite energy condition for $v_.$, one can use the results in \cite{CL1} and build a solution of \eqref{eqannealedbis} satisfying nice properties. The main tool for doing this is Theorem \ref{thmCL} in the Appendix  shown in \cite{CL1} Theorems 4.29, 4.42 and 4.48.

Uniqueness of the solution is not explicitly addressed in \cite{CL1}. This discussion is made in Remark 2.12 of \cite{CatKS}.

We may thus apply all what precedes in our situation.  We gather all this in the following main theorem.
\begin{theorem}\label{thmannealed}
Let $p_t$ be defined by \eqref{eqconvoldens}, and $\hat p_t=p_{\kappa t}$ for some $\kappa \in (0,1)$. Assume the following
\begin{enumerate}
\item[(1)] \quad (i) \; Either $\nabla^2 \ln p_t$ is uniformly bounded on $[0,T^-)$, meaning that $$||\nabla^2 \ln p_t(.)||^2:= \sum_i \, \sup_x  \, \sum_j |\partial^2_{i,j} \ln p_t(x)|^2 \leq a_t^2 \, ,$$ and $\sup_{t \leq T-\varepsilon} a_t=a(\varepsilon) < +\infty$ for all $\varepsilon >0$, 

(ii) or $\nabla \ln p_t$ is uniformly bounded on $[0,T^-)$, i.e. $$\sup_x \sum_i |\partial_i \ln p_t(x)|^2 \leq b^2_t$$ where $\sup_{t \leq T-\varepsilon} b_t=b(\varepsilon) < +\infty$ for all $\varepsilon >0$,
\item[(2)] \quad $\int_0^{T} \, |p'_t|^2  \, dt < + \infty$, $|p'|$ being defined in Theorem \ref{thmambrosio}.
\end{enumerate}
Then there exists a (strongly in case (i), weakly in case (ii)) unique solution of \eqref{eqannealed}, and a (weakly) unique solution of \eqref{eqannealedbis} up to $(T/\kappa)^-$. 

If $Q_X$ and $Q_Y$ denote the distribution on the path space $C^0([0,(T/\kappa)^-),\mathbb R^d)$ of the processes $X_.$ and $Y_.$, $Q_X$ is absolutely continuous w.r.t. $Q_Y$ and $$d_{KL}(Q_X,Q_Y)= \frac 14 \, \int_0^{T/\kappa} \int \,  \, |v_t|^2 \, \hat p_t dx \, dt = \frac{\kappa}{4} \, \int_0^T \, |p'_t|^2 \, dt \, .$$ In addition, the distribution of $X_t$ is exactly $\hat p_t dx$ and we have
\begin{equation}\label{eqfinalKL}
d_{KL}(\tilde \pi,\mathcal L(Y_{T/\kappa})) \leq \frac{\kappa}{4} \, \int_0^T  \, |p'_t|^2 \, dt \, ,
\end{equation}
where $Y_{T/\kappa}$ denotes any weak limit of $Y_{(T/\kappa)-\varepsilon}$ as $\varepsilon \to 0$. If $\lambda_T=1$, $\tilde \pi=\pi$.
\end{theorem}
\begin{proof}
With our assumptions it is immediate that $\hat p_t$ is a weak solution of \eqref{eqFPannealed} i.e. of the weak forward equation in Theorem \ref{thmCL}. Under assumption (1) (i), $\nabla \ln p_t$ is $L$-Lipschitz continuous with $L \leq a(\varepsilon)$ for $t\in [0,T-\varepsilon)$. So is $\nabla \ln \hat p_t$ for $t\in [0,(T-\varepsilon)/\kappa)$. We may thus apply Theorem \ref{thmCL} (2) with $P=Q_Y$. Under assumption (1) (ii), $\nabla \ln p_t + v_t$ is of finite energy and we may apply Theorem \ref{thmCL} with $P$ equal to the Wiener measure with initial distribution $\tilde \nu$. Since $P$ is equivalent (in this situation) to $Q_Y$, $Q_X$ is absolutely continuous w.r.t. $Q_Y$. The density is given by the appropriate Girsanov exponential martingale and the conclusion regarding the Kullback-Leibler divergence between the path-space distributions $Q_X, Q_Y$ follows.

The last inequality then follows immediately since the Kullback-Leibler divergence is non-increasing under measurable transformation and is lower semi continuous w.r.t.  the weak convergence topology. Note that, under our assumptions, it is not clear that $Y_{(T/\kappa)-\varepsilon}$ is converging to some $Y_{T/\kappa}$, this has to be discussed on each example.
\end{proof}

\begin{remark}\label{remequiv}
It is important to remark that $a_t$ in the previous Theorem is exactly what is needed to control the Lipschitz norm of $\nabla \ln p_t$. Of course $$a_t \leq \sup_x \, ||\nabla^2 \ln p_t(x)||_{HS}$$ where $||.||_{HS}$ denotes the Hilbert-Schmidt (or Frobenius) norm of a matrix.

Since $\nabla^2 \ln p_t(x)$ is symmetric, another natural ``pseudo-norm'' is given by $$- \, C_t(x) \, Id \leq \nabla^2 \ln p_t(x) \leq C_t(x) \, Id$$ in the sense of quadratic forms. 

Recall that $C^2_t(x) \leq ||\nabla^2 \ln p_t(x)||_{HS}^2 \leq d \, C^2_t(x)$ hence in case (i), $$a_t \leq  \sqrt d \, C_t \, := \, \sup_x \, \sqrt d \, C_t(x) \, .$$
\hfill $\diamondsuit$
\end{remark}

As is common for interpolations of the form \eqref{eqinterpolationgeneral}, some degeneracy may occur if $\lambda_0=0$ or/and $\lambda_T=1$. Below, we give a result showing how to overcome possible difficulties in these cases.
\begin{corollary}\label{corannealed2}
Define $m_\pi=\int |x| \, \pi(dx)$ and $V_\pi= \int |x|^2 \, \pi(dx)$, and similarly $m_\nu$ and $V_\nu$. Recall that $$\tilde \pi=\mathcal L(\sqrt{\lambda_T} \, X + \sqrt{1-\lambda_T} \, Z) \quad \textrm{ and } \quad \tilde \nu=\mathcal L(\sqrt{\lambda_0} \, X + \sqrt{1-\lambda_0} \, Z ) \, .$$  Then $$W_1(\pi,\tilde \pi) \leq (1-\sqrt{\lambda_T}) m_\pi + \sqrt{1-\lambda_T} \, m_\nu \,\, \textrm{ and } \,\, W_2^2(\pi,\tilde \pi) \leq 2 \, (1-\sqrt{\lambda_T})^2 V_\pi + 2 \, (1-\lambda_T) \, V_\nu \, ,$$ and similar bounds are satisfied for $\nu$ and $\tilde \nu$.

Consequently, $$d_{BL}(\pi,\mathcal L(Y_{T/\kappa})) \leq (1-\sqrt{\lambda_T}) m_\pi + \sqrt{1-\lambda_T} \, m_\nu + \sqrt{\kappa/2} \, \left(\int_0^T \, |p'_t|^2 dt\right)^{\frac 12} \, ,$$ where $d_{BL}$ denotes the bounded Lipschitz distance i.e. $$d_{BL}(\mu,\theta) \, = \, \sup\{\int f d\mu - \int f d\theta \; ; \; ||f||_\infty \leq 1 \textrm{ and } ||\nabla f||_\infty \leq 1\} \, .$$

\end{corollary}
\begin{proof}
The first part is immediate since $(X,\tilde X_T)$ is a coupling of $\pi$ and $\tilde \pi$. Note that if one of $\pi$ or $\nu$ is centered one can skip the factor $2$ for $W_2$. The second part follows from the following facts
\begin{enumerate}
\item $d_{BL} \leq W_1$, since $W_1(\mu,\theta)=\sup\{\int f d\mu - \int f d\theta \, ; \, ||\nabla f||_\infty \leq 1\}$,
\item $d_{BL} \leq d_{TV}$ where $d_{TV}(\mu,\theta)=\sup\{\int f d\mu - \int f d\theta \, ; \, ||f||_\infty \leq 1\}$  denotes the total variation distance,
\item $d_{TV} \leq \sqrt{2 \, d_{KL}}$ which is the celebrated Pinsker inequality.
\end{enumerate}
\end{proof}

\begin{remark}
As we said the above strategy is described in \cite{Guo} but with missing assumptions (see Lemma 1 therein) or references (in particular for \eqref{eqannealedbis}).

Very interesting is that the choice of $v_.$ in Theorem \ref{thmambrosio} is the optimal choice in the set of vector fields satisfying the continuity equation with a finite $\mathbb L^2$ norm. It is shown in \cite{AGS} that this $v_.$ belongs to the $\mathbb L^2$ closure of the gradients. Looking at section 4 in \cite{CL1} one sees that it corresponds to the optimal (minimal entropy) choice of $Q_X$ on the path space among all probability measures with finite relative entropy w.r.t. $Q_Y$ and marginals flow $\hat p_t dx$.
\hfill $\diamondsuit$
\end{remark}

\begin{remark}{\textbf{About the literature.}}

Building a diffusion process with a given flow of time marginals is an old problem. It was stated and solved in \cite{Carlen} for drifted Brownian motion in connection with Nelson stochastics mechanics. The proof in \cite{Carlen} is purely analytic and does not discuss the nature of the path measure (here $Q_X$). The construction in \cite{CL1} is using (delicate) stochastic calculus. Another proof based on large deviations arguments is contained in \cite{CL2,CL3}. The link between Nelson problem and large deviations was pointed out in F\"{o}llmer's lecture notes at Saint Flour \cite{FolSF}. Another related problem is the construction of (Schr\"{o}dinger) bridges where only the initial and the final marginals are given. For this aspect we refer to \cite{FolSF,Leo-survey}. The second reference in particular contains the fundamental relationship between these bridges and entropic transport of measures.
\hfill $\diamondsuit$
\end{remark}

\begin{remark}\label{remlogsobimprove}

If we compare the final bias $d_{KL}(\tilde \pi,\mathcal L(Y_{T/\kappa}))$ with the one obtained with the overdamped Langevin process and time reversal at $T/\kappa$, one can be disappointed. In the annealed case we get something of size $\kappa$ while it is $e^{-T/\kappa}$ for the overdamped case. Of course, if one wants to preserve the fact that the density at time $t$ is explicit, the latter is limited to a gaussian base distribution, while the former does not assume anything on the base distribution, except some regularity for the score function. In addition one does not have to assume that $d_{KL}(\pi,\nu)$ is finite in the annealed case.

One can nevertheless ask whether the final bias can be improved if we assume that $\nu$ satisfies a logarithmic Sobolev inequality, as the gaussian distribution does. We shall come back to this question in Section \ref{seclogsob}, where our main result will be presented.
\hfill $\diamondsuit$
\end{remark}

Sections \ref{secgame}-\ref{secnewpoinc} are devoted to analyzing conditions on $\pi$ and $\nu$ under which Theorem \ref{thmannealed} holds. In particular, we establish how this requirement is connected to controlling the Poincar\'e constant of the conditional distribution. This analysis leads to the study of Poincar\'e inequalities for perturbed measures, where we recall relevant results from \cite{CGperturb} and present new contributions.
\medskip

\section{When functional inequalities enter the game, or not}\label{secgame}

For 
\begin{equation}\label{eqinterpol}
 p_t(x) dx = \frac{\pi(x/\sqrt{\lambda_t})}{\lambda_t^{d/2}}*\frac{\nu(x/\sqrt{1-\lambda_t})}{(1-\lambda_t)^{d/2}} \, ,
\end{equation} 
we thus have to find conditions for assumptions (1) and (2) of Theorem \ref{thmannealed} to be satisfied.
 
We start with (2).
\begin{proposition}\label{propentropambrosio}
Define $\int |x|^2 \, d\pi=V_\pi$ and $\int |x|^2 \, d\nu=V_\nu$. Let $t \mapsto \lambda_t$  be increasing, non-negative and $C^1$ on $[0,T]$. Define $$A_0=\int_0^T \, \frac{|\lambda'_t|^2}{\lambda_t} \, dt \quad ; \quad A_1=\int_0^T \, \frac{|\lambda'_t|^2}{1-\lambda_t} \, dt \, .$$ Then $$\int_0^T |p_t'|^2 \, dt \leq \frac 12 \; (V_\pi \, A_0 \, + \, V_\nu \, A_1) \, .$$
\end{proposition}
\begin{proof}
If $X$ and $Z$ are independent with respective distributions $\pi$ and $\nu$, the pair $$(\sqrt{\lambda_t} X + \sqrt{1-\lambda_t} Z ; \sqrt{\lambda_{t+h}} X + \sqrt{1-\lambda_{t+h}} Z)$$ is a coupling of $p_t$ and $p_{t+h}$ for $t\in(0,T)$ and $h$ small enough. Thus 
\begin{eqnarray*}
W_2^2(p_t,p_{t+h}) &\leq& \mathbb E[\vert(\sqrt{\lambda_{t+h}} - \sqrt{\lambda_t}) X \, + \, (\sqrt{1-\lambda_{t+h}} - \sqrt{1-\lambda_t}) Z\vert^2] \\ &\leq& 2 \mathbb E[(\sqrt{\lambda_{t+h}} - \sqrt{\lambda_t})^2 \vert X\vert^2] \, + \, 2 \, \mathbb E[(\sqrt{1-\lambda_{t+h}} - \sqrt{1-\lambda_t})^2 \vert Z\vert^2]
\end{eqnarray*}
It follows immediately that $$|p'_t|^2 \leq \frac{|\lambda'_t|^2}{2} \, \left(\frac{V_\pi}{\lambda_t}+ \frac{V_\nu}{(1-\lambda_{t})}\right)$$ hence the result.
\end{proof}

\begin{remark}
Notice that the condition on $\lambda_.$ is satisfied for $\lambda_t=2(t/T)^2 \, \mathbf 1_{t \leq T/2} + \, (1-2(1-t/T)^2) \, \mathbf 1_{t > T/2}$ and $\lambda_t=1/2(1+\cos(\pi(1-(t/T)^\alpha)))$ with $\alpha >1/2$.

It is worth noting that the $C^1$ assumption on the schedule can be weakened to $t\mapsto\lambda_t$ being right differentiable.

Moreover, if $X$ or $Z$ is centered, one can replace $1/2$ by $1/4$ in the final estimate.

Finally, the conditions $A_0$ and $A_1$ finite are much weaker than Assumption A6 (or A10) in \cite{Corderoetal}.
\hfill $\diamondsuit$
\end{remark}

We turn to condition (1). For simplicity we shall assume that both $\pi$ and $\nu$ are absolutely continuous w.r.t. the Lebesgue measure so that $$\pi(dx)= e^{-U(x)} \, dx \textrm{ and } \nu(dx) = e^{-W(x)} \, dx \, .$$ We thus have $$p_t(x) = c(t) \, \int \, e^{- U_t(y)} \, e^{- W_t(x-y)} \, dy \, = \, c(t) \, \int \, e^{- U_t(x-y)} \, e^{- W_t(y)} \, dy \, ,$$ with $$c(t)=(\lambda_t \, (1-\lambda_t))^{-d/2} \quad , \quad U_t(z)=U(z/\lambda_t^{1/2}) \quad , \quad W_t(z)=W(z/(1-\lambda_t)^{1/2}) \, . $$ We shall now make formal calculations. They will be justified later. First introduce for all $x$ the conditional probability density 
\begin{equation}\label{eqcondition}
q_t^x(y) = (Z_t^x)^{-1} \, e^{- (U_t(y)+W_t(x-y))} \, , \, Z_t^x = \int e^{- (U_t(y)+W_t(x-y))} dy \,  . 
\end{equation} 
Let $Y_t^x$ denote a random variable with probability density $q_t^x$. One has
\begin{eqnarray*}
\nabla \ln p_t(x)  &=& - \, \int \nabla W_t(x-y) \, q_t^x(y) \, dy = - \, (1-\lambda_t)^{-1/2} \, \mathbb E[\nabla W((x-Y_t^x)/(1-\lambda_t)^{1/2})] \\ &=& - \, \int \nabla U_t(y) \, q_t^x(y) \, dy = \, - \, \lambda_t^{-1/2} \, \mathbb E[\nabla U(Y_t^x/\lambda_t^{1/2})]
\end{eqnarray*}
and
\begin{eqnarray*}
\nabla^2 \ln p_t(x) &=& \, (1-\lambda_t)^{-1} \, \left(- \, \mathbb E[\nabla^2 W((x-Y_t^x)/(1-\lambda_t)^{1/2})] \, + \, \Cov[\nabla W((x-Y_t^x)/(1-\lambda_t)^{1/2})]\right) \\ &=&  \lambda_t^{-1} \, \left( - \, \mathbb E[\nabla^2 U(Y_t^x/\lambda_t^{1/2})] \, + \, \Cov[\nabla U(Y_t^x/\lambda_t^{1/2})]\right) \\ &=& (\lambda_t \, (1-\lambda_t))^{-1/2} \, \Cov[\nabla W((x-Y_t^x)/(1-\lambda_t)^{1/2}),\nabla U(Y_t^x/\lambda_t^{1/2})] \, .
\end{eqnarray*}
\smallskip

\begin{remark}\label{remcommut}
One can of course replace $q_t^x$ by $\bar q_t^x(y)=q_t^x(x-y)=(\bar Z_t^x)^{-1} \, e^{-(U_t(x-y)+W_t(y))}$, $Y_t^x$ by $\bar Y_t^x$ and then exchange $W$ and $U$ (as well as $\lambda_t$ and $1-\lambda_t$) in the preceding formulas. We shall use this remark at several places.
\hfill $\diamondsuit$
\end{remark}

It remains to make some assumptions for the previous calculations to be rigorous.
\begin{assumption}\label{assump}
$e^{-U}$ and $e^{-W}$ are bounded, $V_\pi$ and $V_\nu$ are finite.
\end{assumption}
Under these assumptions, if $|\nabla W|$ (resp. $|\nabla U|$) is bounded, one can differentiate under the integral sign and get the first (resp. second) expression for $\nabla \ln p_t$. 

If $|\nabla^2 W| \, \leq C_W \, Id$, $|\nabla W|$ has at most linear growth, so that thanks to the Assumption \ref{assump} it is integrable w.r.t. $q_t^x(y) \, dy$. We deduce the first expressions for  both $\nabla \ln p_t$ and $\nabla^2 \ln p_t$. If $|\nabla^2 U| \, \leq C_U \, Id$ one obtains the second and the third expressions. 
\medskip

A first result is thus immediate
\begin{theorem}\label{thmbounded}
Assume that $e^{-U}$ and $e^{-W}$ are bounded, that $V_\pi$ and $V_\nu$ are finite, and that $t \mapsto \lambda_t$ is increasing. We will introduce some assumptions

\textbf{(HbW)} \quad $\sup_x |\nabla W(x)| \leq M_W < +\infty$, 

\textbf{(H2bW)} \quad$|\nabla^2 W|\leq C_W Id$, meaning that $- C_W \, Id \leq \nabla^2 W(x) \, \leq C_W \, Id$ in the sense of quadratic forms, for all $x$,

\textbf{(HbU)} \quad $\sup_x |\nabla U(x)| \leq M_U < +\infty$,

\textbf{(H2bU)} \quad $|\nabla^2 U|\leq C_U Id$,  meaning that $- C_U \, Id \leq \nabla^2 U(x) \, \leq C_U \, Id$ in the sense of quadratic forms, for all $x$.
\medskip

\begin{enumerate}
\item[(1)] \quad Assume that (HbW)  is satisfied. Then $$|\nabla \ln p_t| \leq (1-\lambda_t)^{-1/2} \, M_W$$ so that (1) (ii) in Theorem \ref{thmannealed} is satisfied, so that one may apply Theorem \ref{thmannealed} up to time $(T/\kappa)^-$. If in addition (HbU) is also satisfied, one has for all $t$ 
\begin{eqnarray*}
|\nabla \ln p_t| &\leq& \min ((1-\lambda_t)^{-1/2} \, M_W \; ; \; \lambda_t^{-1/2} \, M_U ) \, \\ &\leq& \, \max (\lambda_{1/2}^{-1/2},(1-\lambda_{1/2})^{-1/2}) \, \max(M_W,M_U)
\end{eqnarray*}
 and one may apply Theorem \ref{thmannealed} up to (and including) time $T/\kappa$.
\item[]
\item[(2)] \quad Assume that (HbW) and (H2bW) are satisfied. Then for all $t\in[0,T)$, $$- \, (1-\lambda_t)^{-1} C_W Id \leq \nabla^2 \ln p_t \leq \, (1-\lambda_t)^{-1} \, (C_W+M_W^2) \, Id \, .$$
If in addition (HbU) and (H2bU) are also satisfied, then for all $t \in [0,T]$, $$- \, c \max(C_U,C_W) Id \leq \nabla^2 \ln p_t \leq  c \, \max(C_W+M_W^2,C_U+M_U^2) \, Id \, ,$$ where $c=\max (\lambda_{1/2}^{-1},(1-\lambda_{1/2})^{-1})$. (1) (i) in Theorem \ref{thmannealed} is satisfied and the conclusion of the Theorem is in force.
\end{enumerate}
As a by product one can apply Corollary \ref{corannealed2}. 
\end{theorem}
\begin{remark}\label{remdim1}
As we previously mentioned, (2) is interesting in order to prove strong existence for \eqref{eqannealed}, and also furnishes important useful controls  for discretization schemes. It is worth to note that one can take $\lambda_0=0$ as soon as (HbW)and (H2bW) are satisfied, and also $\lambda_T=1$ if (HbU) and (H2bU) are satisfied.

Also notice that the lower bounds in (2) are still true without assuming (HbW) or (HbU).

Another important point is to get some explicit upper bound for the Lipschitz constant $L_t$ of $\nabla \ln p_t$. What precedes shows that $$L_t \leq   \, \frac{\sqrt d}{1-\lambda_t} \, (C_W+M_W^2) \,\,\, \textrm{ or } \,\,\, L_t \leq \, \frac{\sqrt d }{\min (\lambda_{1/2},(1-\lambda_{1/2}))} \, \max(C_W+M_W^2,C_U+M_U^2) \, .$$
\hfill $\diamondsuit$
\end{remark}

\begin{example}\label{exheavy}
In section 4 of \cite{Corderoetal} the authors consider the case where $\nu$ is given by a multivariate Student distribution, while $\pi$ is a compactly supported perturbation of another multivariate Student distribution.

Recall that the multivariate Student distribution $t(0,\sigma^2 Id,\alpha)$ is given by a density $$q(\alpha,\sigma,d,y) = z^{-1} \, \left(1 + \frac{|y|^2}{\alpha \sigma^2}\right)^{- (\alpha +d)/2}$$ so that one has to assume that $\alpha>2$ for its variance to be finite (the mean being $0$ in this case), i.e. $V_W<+\infty$. One also uses the terminology generalized Cauchy distribution.

If $\nu$ is a $t(0,\sigma^2 Id,\alpha)$, then $$W(y)=\ln z+\frac{\alpha + d}{2} \, \ln\left(1 + \frac{|y|^2}{\alpha \sigma^2}\right) \quad \textrm{ so that } \quad \nabla W(y) = (\alpha + d) \, \frac{y}{\alpha \sigma^2 + |y|^2} \, .$$ In particular $|\nabla W|$ is bounded by $M_W=(\alpha +d)/(2 \, \sigma \sqrt \alpha)$. 

Since $$\partial_{i,j}^2 W(y) = \frac{(\alpha+d)}{\alpha \sigma^2+|y|^2} \; \left(\delta_{i,j} - \frac{2 y_i y_j}{\alpha \sigma^2+|y|^2}\right)$$ it is a simple exercise to check that $$\frac{(\alpha+d)}{\alpha \sigma^2} Id \geq \nabla^2W \geq - \, \frac{(\alpha+d)}{2\alpha \sigma^2} Id \; \textrm{ so that } C_W \leq \frac{(\alpha+d)}{\alpha \sigma^2} \, .$$

One can similarly consider Subbotin (exponential power) distributions $$\nu(dx) = z^{-1} e^{- (1+|x|^2)^{\alpha/2}} \, dx$$ for $0<\alpha<1$. Explicit calculations are left to the reader.
\hfill $\diamondsuit$
\end{example}
\bigskip

Assume that (H2bW) is satisfied. If we do no more assume that $\nabla W$ is bounded, as already remarked, the lower bound in (2) of Theorem \ref{thmbounded} is still true.  For the upper bound we shall follow the idea in \cite{Corderoetal} and use a functional inequality. Indeed on one hand, we have to bound $\nabla^2 W$ and possibly $\nabla^2 U$ from below, and on the other hand we have to bound $\Cov[\nabla W((x-Y_t^x)/(1-\lambda_t)^{1/2})]$ and possibly $\Cov[\nabla U(Y_t^x/\lambda_t^{1/2})]$ from above. This leads to 
\begin{lemma}\label{lemmalower}
Assume that there exists some constant $C_P(t)$ such that for all $x$ the distribution $q^x_t(y)dy$ satisfies a Poincar\'e inequality with constant less than $C_P(t)$, i.e. for all nice function $g$, $$\Var(g(Y_t^x)) \leq C_P(t) \, \mathbb E(|\nabla g|^2(Y_t^x)) \, .$$  
\begin{enumerate}
\item[(1)] \quad Assume that $|\nabla^2 W| \, \leq C_W \, Id$, meaning that $- \, C_W \, Id \leq \nabla^2 W \, \leq C_W \, Id$ in the sense of quadratic forms, so that $||\nabla^2 W||^2 \leq d C^2_W$. Then $$|\nabla^2 \ln p_t| \leq \, \frac{C_W}{1-\lambda_t} \, \left( 1 + \frac{ d \, C_W \, C_P(t)}{1-\lambda_t}\right) \; Id \, .$$
\item[(2)] \quad Assume in addition that $|\nabla^2 U| \, \leq C_U \, Id$. Then $$|\nabla^2 \ln p_t| \leq \, c \, \max (C_W,C_U) \, (1 \, + \, cd \, \max (C_W,C_U) \, C_P(t)) \; Id \, ,$$ where $c$ is as in Theorem \ref{thmbounded}.
\end{enumerate}
Following remark \ref{remcommut} we may replace $q_t^x(y)dy$ by $\bar q_t^x(y)dy$ that satisfies a Poincar\'e inequality with the same constant $C_P(t)$ for all $x$.
\end{lemma}
\begin{remark}\label{remlipdim}
Recall that the Lipschitz constant $L_t$ of $\nabla \ln p_t$ is less than $\sqrt d$ times the constant obtained in the bound for $|\nabla^2 \ln p_t|$. In particular even if $C_P(t)$ is dimension free, the dimension dependence becomes of order at least $d^{3/2}$, except in some special cases discussed in Remark \ref{remconvolgene}. \hfill $\diamondsuit$
\end{remark}

\begin{proof}
It is enough to show that the right hand side is an upper bound for $\nabla^2 \ln p_t$, since it is larger than the (opposite of) the lower bound we already know. The proof is immediate using the straightforward extension of Poincar\'e inequality to multivalued functions yielding (see Lemma B.8 in \cite{Corderoetal}) 
$$\Cov[\nabla W(x-Y_t^x/(1-\lambda_t)^{1/2})]\, \leq \, \sup_{|\xi|=1} \, \Var (\langle \nabla W(x-Y_t^x/(1-\lambda_t)^{1/2}), \xi\rangle)\, Id$$
\begin{eqnarray*}
 &\leq&  C_P(t) \, (1-\lambda_t)^{-1} \, \sup_{|\xi|=1} \, \mathbb E\left[\sum_i \, \big|\sum_j \, \partial^2_{ij} W (x-Y_t^x/(1-\lambda_t)^{1/2}) \, \xi_j\big|^2\right]\, Id\\ &\leq& C_P(t) \, (1-\lambda_t)^{-1} \, \mathbb E[||\nabla^2 W(x-Y_t^x/(1-\lambda_t)^{1/2})||_{HS}^2]\, Id \, .
\end{eqnarray*}
\end{proof}

Everything in what precedes is in a sense satisfactory in order to apply the results in the previous section, except, at a first glance, the existence and overall the control w.r.t. $t$ of $C_P(t)$. In particular, one can take $\lambda_0=0$ in (1) and also $\lambda_T=1$ in (2), provided one can bound $C_P(t)$ in a neighborhood of $0$ or $T$. In the next section we investigate these problems.

\begin{remark}\label{remwhypoinc}
One may ask whether it is really necessary to use a heavy tool like the Poincar\'e inequality in order to get an upper bound for $\Cov[\nabla W((x-Y_t^x)/(1-\lambda_t)^{1/2})]$. Indeed if we assume that $|\nabla^2_{i,j} W|\leq C_W$ for all $(i,j)$, $$|\nabla W(y)|^2 \leq d \, C_W^2 |y|^2 + c$$ for some constant $c$, so that $$\Cov[\nabla W((x-Y_t^x)/(1-\lambda_t)^{1/2})] \leq \frac{d C_W^2 \, \mathbb E[|x-Y_t^x|^2]}{1-\lambda_t} \, + \, c \, .$$ Unfortunately the right hand side a priori depends on $x$. Of course we may exchange the role of $W$ and $U$, but similarly, we will have to control $\mathbb E(|Y_t^x|^2)$ which also depends on $x$.
\hfill $\diamondsuit$
\end{remark}
\medskip

\begin{remark}\label{remconvolgene}
The prototype of $\nu$ entering the framework of the previous lemmata is the gaussian measure $\mathcal N(0, \sigma^2 Id)$, centered for the sake of simplicity. $\nu$ is thus strictly log-concave (meaning that $W$ is strictly convex) i.e. $$C_W \, Id \geq \nabla^2 W \geq D_W \, Id$$ for some $D_W >0$. In this situation the estimates in Lemma \ref{lemmalower} improve as
\begin{equation}\label{eqconvex1}
\nabla^2 \ln p_t \leq  \, \left(- \, \frac{D_W}{1-\lambda_t} \, + \,  \frac{d \, C^2_W \, C_P(t)}{(1-\lambda_t)^2}\right) \; Id \, .
\end{equation}

In the gaussian case one has the stronger 
\begin{equation}\label{eqconvexgauss}
\nabla^2 \ln p_t \leq  \, \frac{1}{\sigma^2(1-\lambda_t)} \, \left( - \, 1 \, + \frac{C_P(t)}{\sigma^2(1-\lambda_t)}\right) \; Id \, .
\end{equation}
Indeed $\Var(\langle \xi, \nabla W(Z)\rangle)=1/\sigma^4\,\Var(\langle \xi,Z\rangle)$ for any unit vector $\xi$, so that one can directly use the second inequality in the proof of Lemma \ref{lemmalower} (without using the final Cauchy-Schwarz inequality). In particular, the dimension dependence disappears in the conclusion of the lemma. 

More generally if $\nu$ is a strictly log-concave product measure with $W(x)=\sum W_i(x_i)$ and  $W_i''\geq (1/\sigma^2_i) >0$, one has 
\begin{equation}\label{eqconvexprod}
\nabla^2 \ln p_t \leq  \, \left(- \, \frac{1}{\max_i \sigma_i^2(1-\lambda_t)} \,  + \, \frac{C_P(t)}{\min_i \sigma_i^4(1-\lambda_t)^2}\right) \; Id \, .
\end{equation}
\hfill $\diamondsuit$
\end{remark}

\section{About the Poincar\'e constant $C_P(t)$}\label{secaboutpoinc}

The literature on the Poincar\'e inequality is almost impossible to master. We shall only try here to understand how to find sufficient conditions for $C_P(t)$ to exist and how it depends, in these situations, on $t$.
\medskip
\subsection{General properties of the Poincar\'e inequality}
First recall some basic facts. The Poincar\'e inequality is written as $$\Var(g(Z)) \leq C_P(Z) \, \mathbb E(|\nabla g|^2(Z))$$ i.e.$$\int g^2(x) \, \mu_Z(dx) \leq C_P(\mu_Z) \, \int |\nabla g|^2(x) \, \mu_Z(dx) + \left(\int g(x) \, \mu_Z(dx)\right)^2$$ where $\mu_Z$ denotes the distribution of the random variable $Z$, $C_P(Z)=C_P(\mu_Z)$ denoting in general the optimal constant in the previous inequality. 
\medskip

The following properties are easily shown:
\begin{equation}\label{eqpoincrecall1}
\textrm{for any $x \in \mathbb R^d$, }\quad  C_P(x+Z)=C_P(Z), 
\end{equation}
\begin{equation}\label{eqpoincrecall2}
\textrm{for any $\lambda \in \mathbb R$, }\quad C_P(\lambda Z)=\lambda^2 \, C_P(Z), 
\end{equation}
if $Z_1$ and $Z_2$ are independent and $\lambda \in [0,1]$
\begin{eqnarray}\label{eqpoincrecall3}
C_P(Z_1,Z_2)&=& \max (C_P(Z_1),C_P(Z_2)) \; , \nonumber \\ C_P(\sqrt \lambda \, Z_1+\sqrt{1-\lambda} \, Z_2)&\leq& \lambda C_P(Z_1)+ (1-\lambda) \, C_P(Z_2) \, . 
\end{eqnarray}
For the latter one can see e.g. \cite{BBN} Proposition 1. One also easily sees that
\begin{equation}\label{eqconv}
\textrm{If $\mu_n$ weakly converges to $\mu$, $C_P(\mu) \leq \liminf C_P(\mu_n)$.}
\end{equation}

\subsection{Convexity and Poincar\'e inequality}
As emphasized by the Bakry–\'Emery curvature-dimension criterion, (strictly) log-concave measures play a special role in the realm of functional inequalities, and in particular Poincar\'e inequality.
\begin{lemma}\label{lemconvex}
Assume that $H=H_1+H_2$ where for all $x$, $\nabla^2 H_1(x) \geq C_1 \, Id$ for some $C_1>0$, and $H_2$ is bounded. Then $$C_P(e^{-H(x)}dx) \leq \frac{e^{\Osc (H_2)}}{C_1} \, ,$$ where $\Osc (h)=\sup h - \inf h$ denotes the oscillation of $h$.
\end{lemma}
The previous Lemma is a simple consequence of Bakry-\'Emery criterion and Holley-Stroock perturbation argument (see \cite{BaGLbook}). The following is known.
\begin{lemma}\label{lemconvexinfty}
If $H$ is $C^2$ and satisfies $\nabla^2 H(x) \geq C \, Id$ for some $C>0$ and all $|x|\geq R$, then one can build $H_1$ and $H_2$ as in Lemma \ref{lemconvex}. Precisely one can obtain $$\nabla^2 H_1(x) \geq (C/2) \, Id,$$ and $$\Osc(H_2) \leq 16 \, R^2 \, L_H \, ,$$ where $L_H$ denotes the Lipschitz constant of $\nabla H$. Recall that $L_H \leq ||\nabla^2 H||$ where $||\nabla^2 H||^2 = \sup_x \sum_i \, \big|\sum_j \, \partial^2_{i,j}H(x)\big|^2 $.
\end{lemma}
There are several proofs of this type of result. A very simple one in is sketched in Annex B of \cite{aline} (with some gaps), we may also mention \cite{Ghomi} where $H_2$ is not only bounded but also compactly supported. The  version stated in the previous Lemma is shown in the supplementary information of \cite{Ma} (see Lemma 1 therein). 

Replacing $H$ by $H+\varepsilon |x|^2$ and letting $\varepsilon$ go to $0$, it is immediately seen that the previous Lemma extends to the case $C=0$. We can thus extend the Lemma to the case $C\leq 0$ as follows.
\begin{lemma}\label{lemconvexinfty2}
Assume that $H$ is $C^2$ and satisfies $\nabla^2 H(x) \geq C \, Id$ for some $C \in \mathbb R$ and all $|x|\geq R$. Then one can build $H_1$ and $H_2$ such that $$\nabla^2 H_1 \geq 0$$ and $$\Osc(H_2) \leq 16 \, R^2 \, (L_H+|C|)$$  where $L_H$ denotes the Lipschitz constant of $\nabla H$, such that for all $y$, $$H(y)= \frac C2 |y|^2 \, + \, H_1(y) \, + \, H_2(y) \, .$$ 
\end{lemma}
\begin{proof}
It is enough to apply the previous Lemma replacing $H(y)$ by $\tilde H(y)=H(y) - (C/2)|y|^2$. The only thing to remark is that $L_{\tilde H} \leq L_H + |C|$. 
\end{proof}

We may thus state a first result 
\begin{theorem}\label{thmpoincmut}
Let $p_t(x)dx$ be given by \eqref{eqinterpol} and $q_t^x(y)$ (resp. $\bar q_t^x(y)$) be given by \eqref{eqcondition} (resp. remark \ref{remcommut}). Assume that $t \mapsto \lambda_t$ is increasing. Let $R\geq 0$. 

Assume that for all $|x|\geq R$, $\nabla^2 W(x) \geq D^R_W \, Id$ and $\nabla^2 U(x) \geq D^R_U \, Id$ for some $D_W^R$ and $D^R_U$ in $\mathbb R$. Define $$c^R(t):= \frac{D^R_W}{1-\lambda_t}+\frac{D^R_U}{\lambda_t} \, ,$$ and $t^R_W=\inf \{t \; ; \; c^R(t)>0 \}$ and  $t^R_U=\sup \{t \; ; \; c^R(t)>0 \}$.

Then, for all $t$ such that $t^R_U>t>t^R_W$, both $q_t^x(y) dy$ and $\bar q_t^x(y) dy$ satisfy a Poincar\'e inequality with $$C_P(q_t^x(y) dy) \leq \frac{1}{c^R(t)} \, e^{16 \, R^2 \, (L_W+ |D_W^R|+L_U+|D_U^R|)} \, ,$$ 
where $L_W$ and $L_U$ denote the corresponding Lipschitz constants.
\end{theorem}
\begin{proof}
According to the previous Lemma, for all $z$, we may decompose $$W(z)=\frac{D_W^R}{2} \, |z|^2 + W_1(z) + W_2(z)$$ where $W_1$ is convex and $\Osc(W_2) \leq 16 \, R^2 \, (L_W +|D_W^R|)$. We have a similar decomposition for $U$. We may use these decompositions with $z=x-y$ for $W$ and $z=y$ for $U$.

It follows, for any $y$ and any $x$,  $$W_t(x-y)+U_t(y) = (1-\lambda_t)^{-1} D_W^R \, \frac{|x-y|^2}{2} + \lambda_t^{-1} \, D_U^R \, \frac{|y|^2}{2} \, + \, \theta^x_t(y) \, + \, \eta^x_t(y)$$ where for all $x$, $\nabla^2 \theta^x_t(.) \geq 0$ and $\Osc(\eta^x_t) \leq 16R^2 (L_W + |D_W^ R| + L_U + |D_U^R|)$. It thus holds $$W_t(x-y)+U_t(y) = A^x_t(y) + \, \eta^x_t(y)$$ where for all $x$, $\nabla^2 A_t^x \geq c^R(t) \, Id$. 

The result for $q_t^x(y) dy$ follows from Lemma \ref{lemconvex}. The case of $\bar q_t^x(y) dy$ is similar.
\end{proof}

\begin{remark}\label{remconvex}
The decomposition with an explicit quadratic term is crucial for the proof. It yields a loss on the perturbation constant but allows a simple result.

If both $D^R_W$ and $D_U^R$ are non-negative, one of them being positive, one has $t_W^R=0$ and $t_U^R=T$. If $D^R_W>0$ but $D_U^R \leq 0$, $t_U^R= T$ but $t_W^R=0$ only if $\lambda_0 \geq |D_U^R|/(D^R_W+|D_U^R|)$. Note that $W_1(\tilde \nu,\nu) \leq \sqrt \lambda_0 \, m_\pi + (1-\sqrt{1-\lambda_0}) \, m_\nu$, so that one can control the initial bias. 
\hfill $\diamondsuit$
\end{remark}
\begin{corollary}\label{corconvex}
Assume that $e^{-W}$ and $e^{-U}$ are bounded, that $V_\pi$ and $V_\nu$ (see Proposition \ref{propentropambrosio}) are finite and that $t \mapsto \lambda_t$ is increasing. 

Assume in addition that $D_U \, Id \leq \nabla^2 U \leq C_U \, Id$ and that $D_W^R Id \leq \nabla^2 W(x) \leq C_W^R Id$ for some $R\geq 0$, all $|x|\geq R$ and some $C_W^R,D_W^R>0$ (which imply that $|\nabla^2 W|\leq C'_W Id$). 

Finally assume that $\lambda_0(D^R_W-D_U)+D_U>0$.

Then \eqref{eqfinalKL}, Theorem \ref{thmannealed}, applies with $\tilde \pi=\pi$, the right hand side being given by Proposition \ref{propentropambrosio}.
\end{corollary}
\begin{remark}\label{remexchange}
In the previous Remark and Corollary we pushed forward the properties of $W$ since it is chosen by the user, while the data have to fit $U$. One can of course exchange the roles of $W$ and $U$ (exchanging $\lambda_t$ and $1-\lambda_t$ too) in the previous statements.
\hfill $\diamondsuit$
\end{remark}
This result contains the main result in \cite{Corderoetal}. If one wants $\lambda_0=0$, one has to assume convexity at infinity of $U$, as in \cite{Corderoetal} Assumption A.4. If $\lambda_0>0$, from a practical point of view one has first to estimate $p_0(x)$. Since one may choose $W$, one can choose a large $D_W^R$ so that $\lambda_0$ is small, possibly equal to the step size of some discretization scheme. The counterpart is that the bounds for $\nabla^2 \ln p_t$ deteriorate. 
\bigskip

\subsection{Perturbation of Poincar\'e inequality}
The previous Theorem is based on convexity results for $W+U$. We will now adopt a perturbative point of view: starting from $\nu_t^x(y)dy := e^{-W((x-y)/\sqrt{1-\lambda_t})} dy$, look at $q_t^x$ as a perturbation of $\nu_t^x$. This approach is studied in \cite{CGperturb}. We shall give below the most tractable existing results we know in this direction. We will come back to the perturbation point of view in the next section.
\begin{theorem}\label{thmperturb}
Let $p_t(x)dx$ be given by \eqref{eqinterpol} and $q_t^x(y)$ be given by \eqref{eqcondition}. Assume that $t \mapsto \lambda_t$ is increasing. 
\begin{enumerate}
\item[(1)] \quad Assume that $\nabla^2 W \geq D_W Id$ with $D_W>0$ and that $|\nabla U|\leq M_U$. Then $$C_P(q_t^x(y)dy) \, \leq \, \frac{2 \, (1-\lambda_t)}{D_W } \, e^{4\sqrt{2d/\pi} \, \frac{ M_U^2 \, (1- \lambda_t)}{\lambda_t  \, D_W}} \, .$$ (Due to Miclo, see lemma 2.1 in \cite{BGMZ} recalled in Theorem 1.3 (1) of \cite{CGperturb} with a typo.)
\item[]
\item[(2)] \quad Under the same assumptions, $$C_P(q_t^x(y)dy) \, \leq  \, 2\left(\frac{M_U \, (1- \lambda_t)}{\sqrt{\lambda_t}  \, D_W}+ \sqrt{ \frac{2 \, (1- \lambda_t)}{D_W}}\right)^2 \; e^{\frac{M_U^2 \, (1-\lambda_t)}{2 \lambda_t  \, D_W}} \, .$$ (Example 3 section 7.1 in \cite{CGsemin}, recalled in Theorem 1.3 (2) of \cite{CGperturb}.)
\item[]
\item[(3)] \quad If $\nabla^2 W(x) \geq D_W^R$ for some $D_W^R>0$ and all $|x|\geq R$, and $|\nabla U|\leq M_U$, $$C_P(q_t^x(y)dy) \, \leq \, \frac{4 \, (1-\lambda_t)}{D_W^R} \, e^{16 \sqrt{2d/\pi} \, \frac{M_U^2 \, \, (1- \lambda_t) }{\lambda_t \, D_W^R}} \, e^{16 R^2 \, L_W} \, ,$$ and $$C_P(q_t^x(y)dy) \, \leq \, 8 \, \left(\frac{M_U \, (1- \lambda_t) }{\sqrt{\lambda_t} \, D_W^R}+ \sqrt{ \frac{1-\lambda_t}{D_W^R}}\right)^2 \; e^{\frac{ M_U^2 \, (1-\lambda_t)}{\lambda_t  \, D_W^R}} \, e^{16 R^2 \, L_W} \, .$$
\item[]
\item[(4)] \quad Assume that for some $\varepsilon >0$, $$s:=\frac{(1+\varepsilon) (1-\lambda_t) }{4 \, \lambda_t} \, C_P(\nu) \, M_U^2 \, < 1 \, ,$$ then $$C_P(q_t^x(y)dy) \, \leq \, \frac{1+(1/\varepsilon)}{1-s} \, (1-\lambda_t) \, C_P(\nu) \, .$$ (See Theorem 2.1 in \cite{CGperturb} and recall that $\nu(dy)=e^{-W(y)} dy$.)
\end{enumerate}
\end{theorem}
\begin{corollary}\label{corperturb}
Under the basic assumptions of Corollary \ref{corconvex}, the assumptions in any of the items of the previous Theorem and the additional assumptions $|\nabla^2 W| \leq C_W Id$ (if not already implied)  and $|\nabla^2 U| \leq C_U Id$, \eqref{eqfinalKL}, Theorem \ref{thmannealed}, applies with $\tilde \pi=\pi$, the right hand side being given by Proposition \ref{propentropambrosio}.
\end{corollary}
\begin{remark}\label{remexchange2}
Again, in the previous Theorem and Corollary we pushed forward the properties of $W$ since it is chosen by the user, while the data have to fit $U$. 
\hfill $\diamondsuit$
\end{remark}

To get (3) from (1) and (2) in the Theorem, simply use Lemma 4.2 and Holley-Stroock argument, for (4) use \eqref{eqpoincrecall1} and \eqref{eqpoincrecall2} saying that $C_P(\nu_t^x(y)dy)=(1-\lambda_t)C_P(e^{-W(y)} dy)$. 

Even if the value of the Poincar\'e constant can be desperately big if $\lambda_0$ is close to $0$, the first three items of the previous Theorem complete Theorem \ref{thmpoincmut} (except if $U$ is strictly convex at infinity too) for Lipschitz $U$. In addition there is no restriction on $t$.
\medskip

(4) does not require any convexity assumption but it is limited to a small range of $t$'s. Indeed the condition $s<1$ amounts to 
\begin{equation}\label{eqlambmin}
\lambda_t > \frac{(1+\varepsilon)C_P(\nu)M_U^2}{4+ (1+\varepsilon)C_P(\nu)M_U^2} \, :=\lambda_{min} \, .
\end{equation}

\begin{remark}\label{remconv+lip}
The previous Theorem completes the range of potential examples, including for instance the case of a gaussian $\nu$ with covariance $\sigma^2 \, Id$ and an heavy tailed $\pi$. Recall that $C_P(\nu)=\sigma^2$.

In this case, in order to get an upper bound for $\nabla^2 \ln p_t$ one has to find an upper bound for $\Var(\langle \xi,(x-Y_t^x)\rangle)=\Var(\langle \xi,Y_t^x\rangle)$ for any unit vector $\xi$. If $\pi$ is compactly supported in the (euclidean) ball $B(0,R)$, the latter is bounded by $\lambda_t R^2$. Notice that in this case we can drop the regularity assumption for $\pi$. We thus have
\begin{proposition}\label{propgausscomp}
Let $\nu$ be a $\mathcal N(0,\sigma^2 Id)$ and $\pi$ be compactly supported in the euclidean ball $B(0,R)$. Then for all $t$, $$- \, \frac{1}{\sigma^2 \, (1-\lambda_t)} \, Id \, \leq \nabla^2 \ln p_t \, \leq \, - \, \frac{\sigma^2 \, (1-\lambda_t) - \lambda_t R^2}{\sigma^4 \, (1-\lambda_t)^2} \, Id \, .$$
\end{proposition}
This result is not new and is (up to the presentation) contained in subsection 2.1 of \cite{BGMZ}.
The previous proposition can be extended to the case where $\pi$ is given by the convolution of a gaussian with a compactly supported distribution.
\begin{proposition}\label{propgausscompgauss}
Let $\nu$ be a $\mathcal N(0,\sigma^2 Id)$ and let $\pi$ be the convolution of $\mathcal N(0,\tau^2 Id)$ with a distribution supported in the euclidean ball $B(0,R)$. Then for all $t$, $$- \, \frac{1}{\alpha_t^2} \, Id \, \leq \nabla^2 \ln p_t \, \leq \, - \, \frac{\alpha_t^2 - \lambda_t R^2}{\alpha_t^4} \, Id \, ,$$
where $\alpha_t^2 =\sigma^2 \, (1-\lambda_t) + \tau^2\lambda_t$.
\end{proposition}
\begin{proof}
The result follows from the observation that $X_t = \sqrt{\lambda_t} X + \sqrt{1-\lambda_t} \sigma Z\overset{d}{=} \sqrt{\lambda_t} U + \sqrt{(1-\lambda_t)\sigma^2 + \lambda_t\tau^2} Z'$  where $X \sim\pi$, $Z, Z'\sim\mathcal{N}(0, I)$, $U$ is compactly supported, and by applying Proposition~\ref{propgausscomp}.
\end{proof}
\hfill $\diamondsuit$
\end{remark}

\begin{remark}\label{remdim}
An important point in Statistics is the dimension dependence of all constants. In the strictly convex case estimates are dimension free, while in the convex at infinity case $||\nabla W||$ grows linearly with $d$. Notice however that $V_\pi$ and $V_\nu$ also grow linearly. 

In case (4) if we choose $\nu$ as a product measure, the Poincar\'e constant is dimension free. If it is log-concave (but not strictly log-concave) it is known that the Poincar\'e constant at most grows as $\ln d$, see \cite{KlarKLS}. Log-concave distributions are now playing an important role in Statistics (see e.g. \cite{SW,Chewi}).
\hfill $\diamondsuit$
\end{remark}
\medskip

\section{New perturbation results for the Poincar\'e constant}\label{secnewpoinc}

The perturbative results we used for Theorem \ref{thmperturb} are based on various techniques for proving a Poincar\'e inequality: (1) relates a Lipschitz perturbation to a bounded one, (2) uses the equivalence between Poincar\'e inequality and exponential $\mathbb L^2$ convergence to equilibrium for the semi-group, the latter being obtained via the reflection coupling method introduced in \cite{EbeCRAS,Ebe}, (4) is using a direct rough approach.
\medskip

Another method was introduced in \cite{BBCG, BCG}: the use of (Foster)-Lyapunov functions. Foster-Lyapunov functions are a central tool in the MCMC community, and their use for exponential stabilization of diffusion processes was pointed out in \cite{DMT}. Recall two main results for the Poincar\'e inequality.
\begin{theorem}\label{thmlyappoinc}
Let $\mu(dx) = e^{-V(x)} dx$, for some smooth $V$, be a probability measure. Denote $L_V:=\Delta - \nabla V.\nabla$.

Assume that there exists a $C^2$ function $F$ such that $F\geq 1$, some $R>0$, constants $b\geq 0$ and   $\theta>0$ such that $$L_VF(x) \leq - \theta F(x) + b \, \mathbf 1_{|x|\leq R} \, ,$$ such a $F$ is called a Lyapunov function (for $L_V$).

Then $\mu$ satisfies a Poincar\'e inequality with constant $$C_P(\mu) \leq \frac{1}{\theta} \, (1 + C_P(\mu_R))$$ where $\mu_R(dx)= \frac{e^{-V(x)} \, \mathbf 1_{|x|\leq R}}{\mu(B(0,R)} \, dx $ is the (normalized) restriction of $\mu$ to the euclidean ball $B(0,R)$.

If in addition $\frac{\partial F}{\partial n} \leq 0$ on $S_R:=\{|x|=R\}$, where $n$ denotes the inward normal to $S_R$, then $$C_P(\mu) \leq \frac{1}{\theta} + C_P(\mu_R) \, .$$

For $d\geq 2$, $$C_P(\mu_R) \leq \, \frac{d+2}{d(d-1)} \, R^2 \, e^{\Osc_{B(0,R)} V} \, ,$$ while for $d=1$ the pre-factor has to be replaced by $4/\pi^2$.

\end{theorem}
The first part is Theorem 1.4 in \cite{BBCG} while the second one is shown in \cite{CGZ} Remark 3.3. The control of $C_P(\mu_R)$ is a consequence of Holley-Stroock perturbation result and the known upper bound (asymptotically sharp as $d \to + \infty$) for the Poincar\'e constant of the uniform measure on $B(0,R)$ (see \cite{BJM} Corollary 4.1 or \cite{CGWradial} Example 5.5). The case $d=1$ is well known.

Actually, as shown in Theorem 2.3 of \cite{CGZ}, there is an equivalence between the existence of a Lyapunov function and the finiteness of the Poincar\'e constant. One can also control $\theta$ in the definition of a Lyapunov function by the Poincar\'e constant, see Proposition 3.1 in \cite{CGZ} and Theorem 2.1 in \cite{CGhit}.

As a direct consequence, if $V=W+U$ is such that $F$ is a $L_W$ Lyapunov function, it is also a $L_V$ Lyapunov function provided 
\begin{equation}\label{eqlyappert}
(-\nabla U.\nabla F)(x) \leq \theta' \, F(x) \quad \textrm{ for } \, |x|>R' \; \textrm{ for some } R'>R \, \textrm{ and } \, \theta'<\theta \, .
\end{equation}
Of course in the framework of the previous sections, one has to consider $V(y) = U(x-y)+W(y)$ (recall that the Poincar\'e constant is unchanged by translation) and get a result independent of $x$. Clearly the only tractable situation is the one where $U$ is Lipschitz continuous with $|\nabla U|\leq M_U$ and such that $$ M_U \, |\nabla F(x)| \leq \theta' \, F(x) \quad \textrm{ for $|x|$ large, with $\theta' < \theta$} \, .$$ 

We shall focus on a class of examples, namely potentials $W$ satisfying some ``almost convex'' property
\begin{equation}\label{eqpseudoconv}
\textrm{there exist $\alpha>0$, $\beta\geq 1$ and $R\geq 0$ such that for $|x|\geq R \; , \; \langle x,\nabla W(x)\rangle \geq \alpha |x|^\beta$.}
\end{equation}
Prototypical examples are Subbotin distributions, i.e. $W(x)=|x|^{\beta}$ for $\beta \geq 1$, so that $\alpha=\beta$ and $R=0$. 
\begin{remark}\label{remconvbbcg}
More generally, as shown in \cite{BBCG} Lemma 2.2, if $W$ is convex ($\nabla^2 W \geq 0$), it satisfies \eqref{eqpseudoconv} for $\beta=1$. 

We will revisit the mentioned proof in a quantitative perspective, assuming for simplicity that $W(0)=\min_y W(y)$ which is not a restriction. First let $M=\sup_{|y|\leq 1} |\nabla W(y)|$ so that $W(y)-W(0)\leq M$ for $|y|\leq 1$. Denote $A_M=\{W(y)-W(0)\leq M\}$, so that $B(0,1) \subseteq A_M$. Since $e^{-W}$ is a probability density, it holds $vol(A_M) \leq e^{|W(0)|+M}$ where $vol$ denotes Lebesgue volume. For $x \in A_M$, since $A_M$ is convex, the cone with basis $B_{d-1}(0,1)$ (the unit ball in dimension $d-1$) and vertex $x$ is a subset of $A_M$. The volume of this cone is thus less than the one of $A_M$, i.e. $$\frac{vol(B_{d-1}(0,1)) \, |x|}{d} \leq e^{|W(0)|+M} \; \textrm{ i.e. } \; |x| \leq \frac{d \, \Gamma((d+1)/2) \, e^{|W(0)|+M}}{\pi^{(d-1)/2}}:= R_W  \, .$$ If $|u|= R_W+1$, $u$ does not belong to $A_M$, so that using that $t \mapsto \frac 1t \, (W(tu)-W(0))$ is non-decreasing since $W$ is convex we obtain for $|x| \geq R_W+1$ that $W(x)-W(0)\geq \frac{M}{1+R_W} \, |x|$. Since $\langle x,\nabla W(x)\rangle \geq W(x)-W(0)$ thanks to convexity again, we thus have obtained that \eqref{eqpseudoconv} is satisfied with $\beta=1$ and $\alpha=M/(1+R_W)$. 

Notice that if the previous proof is elementary, the dimension dependence $\alpha \sim (d/c)^{-d}$ is a disaster  for large $d$. Using more sophisticated results one can greatly improve upon the dimension dependence, the price to pay being that quantities are not explicit.

First of all, according to \cite{KlarKLS}, there exists a universal constant $C_{Klar}$ such that if $\nu=e^{-W(x)} dx$ is a log-concave probability distribution, $$C_P(\nu) \leq C_{Klar}  \; \sigma^2(\nu) \;  (1 + \ln d)$$ where $\sigma^2(\nu)$ is the largest eigenvalue of the covariance matrix of $\nu$. It is not easy (but probably not impossible) to trace $C_{Klar}$ in the various papers needed to get the previous result. It is conjectured (Kannan-Lovasz-Simonovits, for short KLS conjecture) that the bound is actually dimension free i.e.  $$C_P(\nu) \leq C_{KLS}  \; \sigma^2(\nu) $$ for some universal constant. Klartag's result is the most recent step in the direction of this conjecture. In some specific situations, KLS conjecture is true: of course if $\nu$ is a product measure using the tensorization property of the Poincar\'e constant or if $\nu$ is radial (or spherically symmetric and $d\geq 2$) where $C_{KLS} \leq 2$ according to Theorem 1.2 in \cite{BJM} (actually one can choose $d/(d-1)$ in dimension $d\geq 2$). It is also true for some families of explicit measures (like the uniform measures on $l^p$ balls for instance).

It thus immediately follows from Theorem 2.1 in \cite{CGhit} that, for any $R>0$,  there exists a $L_W$ Lyapunov function $F_R$ with $$\theta_R = \nu(B(0,R)) \, \min\left(\frac 18 \, , \, \frac{1}{4 C_{Klar} \sigma^2(\nu) \, (1+\ln d)}\right) \, .$$ This time however $F_R$ is not explicit at all, since it is obtained via Lax-Milgram theorem.
\hfill $\diamondsuit$
\end{remark}

If $W$ satisfies \eqref{eqpseudoconv}, $F(x)=e^{\gamma |x|}$ is a $L_W$ Lyapunov function. Indeed 
\begin{eqnarray*}
L_W F &=& \gamma \, \left(\frac{d-1}{|x|} + \gamma \, - \langle \frac{x}{|x|},\nabla W\rangle\right) \, F(x) \\ &\leq& \gamma \, \left(\frac{d-1}{|x|} + \gamma \, - \, \alpha \, |x|^{\beta-1}\right) \, F(x) \, ,
\end{eqnarray*}
for $|x|\geq R$. Choosing $R'\geq R$ large enough and $\gamma$ small enough for 
$$\theta = \gamma \, \left(\alpha (R')^{\beta-1} \, - \gamma - \frac{d-1}{R'}\right) > 0$$ we thus have built a Lyapunov function for $|x|\geq R'$. Notice that $\partial F/\partial n$ is non-positive on any sphere $S_A$. 
\medskip

With all these points in mind we may state several perturbation results
\begin{theorem}\label{thmperturblyap}
Let $\mu(dx)=e^{-V(x)} dx=e^{-(U+W)(x)} dx$, for smooth $U$ and $W$ be a probability measure. Assume that 
\begin{enumerate}
\item[(1)] \quad $W$ satisfies the ``quasi-convex'' condition \eqref{eqpseudoconv} for some $\alpha_W>0$, $\beta_W\geq 1$ and $R \geq 0$, 
\item[(2)] \quad  $U$ satisfies $|\langle x,\nabla U(x)\rangle| \leq \kappa_U \, |x|^{\beta_U}$ for some $\kappa_U \geq 0$, $\beta_U \leq \beta_W$ and $|x|\geq R$. 
\end{enumerate}
Then, if for some $R'\geq R$ one can find $\gamma>0$ such that $$c(W,U):= \, \alpha_W \, (R')^{\beta_W-1} \, - \, \frac{d-1}{R'} - \gamma \, - \, \kappa_U (R')^{\beta_U-1} \,   > \, 0 \, ,$$ $\mu$ satisfies a Poincar\'e inequality with $$C_P(\mu) \leq \frac{1}{\gamma \, c(W,U)} + C_P(\mu_{R'}) \, ,$$ where a bound for $C_P(\mu_{R'})$ is given in Theorem \ref{thmlyappoinc}.

If $\beta_U<\beta_W$ such a $\gamma$ always exists, while if $\beta_U=\beta_W$ it exists provided $\alpha_W > \kappa_U$.
\end{theorem}
Of course we can add an optimization step by increasing $R'$, hence decreasing $1/\gamma c$ but increasing $C_P(\mu_{R'})$. This is left to the courageous reader.

Remark that, for $c(W,U)$ to be positive we need $R'>((d-1)/\alpha_W)^{1/\beta_W}$, so that in Theorem \ref{thmlyappoinc} the bound for $C_P(\mu_{R'})$ is of order $d^{\frac{2-\beta_W}{\beta_W}}$ times the exponential of the Oscillation of $V$. Note that the larger $\beta$ is, the larger the Oscillation of $W$ is. If $W$ is convex we can (in general) improve this dependence.
\begin{theorem}\label{thmperturblyapconv}
In addition to the assumptions of Theorem \ref{thmperturblyap} assume that $W$ is convex. Denote $\nu(dx) = e^{-W(x)} dx$, assumed to be a probability measure as well, and for simplicity a centered probability measure. Then $$C_P(\mu_{R'}) \, \leq \,  \frac{C_{Klar} \, \sigma^2(\nu) \, (1+\ln d)}{\nu(B(0,R'))} \, e^{\Osc_{B(0,R')} U} \, .$$ If in addition $W$ is radially symmetric we may replace $C_{Klar} \, (1+\ln d)$ by $2$.
\end{theorem} 
\begin{proof}
Since $\nu$ is log-concave and $B(0,R')$ is convex, $\nu_{R'}$ is still log-concave. As we recalled in Remark \ref{remconvbbcg}, $C_P(\nu_{R'}) \leq C_{Klar} \sigma^2(\nu_{R'}) \, (1+\ln d)$ and we may replace $C_{Klar} (1+\ln d)$ by $2$ in the radially symmetric case which is preserved by restricting $\nu$ to a ball. Now, if $\xi$ is a unit vector $$\Cov_{\nu_{R'}}(\langle \xi,X\rangle) \leq \frac{\mathbb E_\nu(\langle \xi,X\rangle^2)}{\nu(B(0,R'))} = \frac{\Cov_\nu(\langle \xi,X\rangle^2)}{\nu(B(0,R'))} \leq \frac{\sigma^2(\nu)}{\nu(B(0,R'))} \, .$$ It remains to apply Holley-Stroock bounded perturbation to get the result.
\end{proof}
In a sense the previous results complete the perturbation result by Miclo we recalled and used in (1) of Theorem \ref{thmperturb}, for log-concave but non strictly log-concave $\nu$'s. It turns out however that in the uniformly convex case the method also furnishes new estimates.
\begin{theorem}\label{thmperturblyapstrictconv}
Assume that $\nabla^2 W \geq D_W \, Id$ for some $D_W>0$ and that $|\nabla U|\leq M_U$. Then $$C_P(\mu) \leq \frac{1}{D_W} \, \left(1+e^{2 M_U R/\sqrt{D_W}}\right)$$ with $$R=\frac 12 \, \left((2+(M_U/\sqrt{D_W})) + \sqrt{(2+(M_U/\sqrt{D_W}))^2+4(d-1)}\right) \, .$$
\end{theorem}
\begin{proof}
Assume first that $D_W=1$. We may apply Theorem \ref{thmperturblyap} with $R=0$, $\beta_W=2$, $\alpha_W=D_W$, $\beta_U=1$ and $\kappa_U=M_U$. Choose arbitrarily $\gamma=1 \leq c(W,U)$. It yields $$R' \geq \frac 12 \, \left((2+M_U) + \sqrt{(2+M_U)^2+4(d-1)}\right) \, .$$ This time $\nu_{R'}$ is strictly log-concave so that $C_P(\nu_{R'})\leq 1$. Hence $$C_P(\mu) \leq 1 + e^{\Osc_{B(0,R')} U} \leq 1 + e^{2M_U R'} \, .$$ It remains to use the dilation property of the Poincar\'e constant using the change of variable $x \mapsto x/\sqrt{D_W}$ and replace $M_U$ by $M_U/\sqrt{D_W}$.
\end{proof}
\begin{remark}\label{remcompare}
This result is better than the one by Miclo $$C_P(\mu) \leq \frac{2}{D_W} \, e^{4 \sqrt{2d/\pi} \, M_U^2/D_W}$$ for large values of $M_U/\sqrt{D_W}$ or large dimension $d$. It is however worse than the one in \cite{CGsemin} we recalled $$C_P(\mu) \leq 2\left(\frac{M_U}{D_W} + \sqrt{\frac{2}{D_W}}\right)^2 \, e^{M_U^2/2 D_W} \, ,$$ when the dimension is much larger than $M_U/\sqrt{D_W}$. In general they are not comparable. Notice that the previous result is the only one where the dimension $d$ does not appear in the controls.
\hfill $\diamondsuit$
\end{remark}

Coming back to the situation of the preceding section, we may state
\begin{corollary}\label{cornewpoinc}
Assume that $|\nabla U|\leq M_U$. Then Theorem \ref{thmperturblyap}, Theorem \ref{thmperturblyapconv} and Theorem \ref{thmperturblyapstrictconv} apply to $C_P(\bar q_t^x(y) dy)$ with $\beta_U=1$, replacing $\alpha_W$ ($=D_W$ if $\beta_W=2$), $R$, $\kappa_U=M_U$, $\sigma^2(\nu)$ respectively by $\alpha_W/(1-\lambda_t)^{\beta_W/2} \, $, $R \sqrt{1-\lambda_t}$, $M_U/\sqrt{\lambda_t}$ and $(1-\lambda_t) \, \sigma^2(\nu)$.
\end{corollary}
Indeed $|\langle y,\nabla_y U(x-y)\rangle| \leq M_U |y|$ for all $x$, so that all the previous results apply, uniformly in $x$ to $(\bar Z_t^x)^{-1} \, e^{-(W_t(y)+U_t(x-y))} dy$.
\medskip

All our results now almost cover the full spectrum of typical noise distributions, heavy tailed distributions whose Hamiltonians have bounded gradients, strictly log-concave distributions as the gaussian one and now intermediate distributions. In particular, we may apply one of our results to any Subbotin distribution with any power $0<\alpha\leq 2$. 
\medskip

\section{Using logarithmic Sobolev inequalities}\label{seclogsob}

A measure $\mu$ satisfies a logarithmic Sobolev inequality if the following holds for all $f$ nice enough,
\begin{equation*}
    \Ent_\mu(f^2):=\int f^2\ln\left(\frac{f^2}{\int f^2 d\mu}\right)d\mu \, \leq \, 2C_{LS}(\mu) \, \int |\nabla f|^2 \, d\mu \, = \, 2C_{LS}(\mu) \, \mathbb{E}[\vert \nabla f\vert^2],
\end{equation*}
where $C_{LS}(\mu)>0$ denotes the optimal constant in the previous inequality. We shall use the equivalent formulation $$\int \rho \, \ln(\rho) \, d\mu \, \leq \, (C_{LS}(\mu)/2) \, \int \frac{|\nabla \rho|^2}{\rho}  \, d\mu $$ if $\rho$ is a nice density of Probability. One recognizes the Fisher information of $\rho$ in the right hand side. If the random variable $Z$ satisfies $\mathcal L(Z)=\mu$, we will denote $C_{LS}(Z)$ instead of $ C_{LS}(\mu)$ for convenience.
\medskip

Similar to the Poincar\'e inequality the following properties hold:
\begin{equation}\label{eqlogsobrecall1}
\textrm{for any $x \in \mathbb R^d$, }\quad  C_{LS}(x+Z)=C_{LS}(Z), 
\end{equation}
\begin{equation}\label{eqlogsobrecall2}
\textrm{for any $\lambda \in \mathbb R$, }\quad C_{LS}(\lambda Z)=\lambda^2 \, C_{LS}(Z), 
\end{equation}
if $Z_1$ and $Z_2$ are independent and $\lambda \in [0,1]$
\begin{eqnarray}\label{eqlogsobrecall3}
C_{LS}(Z_1,Z_2)&=& \max (C_{LS}(Z_1),C_{LS}(Z_2)) \; , \nonumber \\ C_{LS}(\sqrt \lambda \, Z_1+\sqrt{1-\lambda} \, Z_2)&\leq& \lambda C_{LS}(Z_1)+ (1-\lambda) \, C_{LS}(Z_2) \, . 
\end{eqnarray}
For the latter one can see e.g. \cite{CGsemin} Proposition 18.  
One also easily sees that
\begin{equation}\label{eqconvlogsob}
\textrm{If $\mu_n$ weakly converges to $\mu$, $C_{LS}(\mu) \leq \liminf C_{LS}(\mu_n)$.}
\end{equation}
\smallskip

We shall try to use here a log-Sobolev inequality in the spirit of the semi-group proof of exponential decay of relative entropy (Kullback-Leibler distance). Our aim is to improve upon the bound 
\begin{equation*}
d_{KL}(\tilde \pi,\mathcal L(Y_{T/\kappa})) \leq \frac{\kappa}{4} \, \int_0^T \, |p'_t|^2 \, dt
\end{equation*}
given in \eqref{eqfinalKL}, Theorem \ref{thmannealed}. Indeed this bound is far from being sharp, as it is a consequence of the decay of the Kullback-Leibler distance by any measurable push forward.
The following theorem presents the main result of this section.

\begin{theorem}\label{thmLSimprov}
    Let $\hat p_t$ denote the density at time $t$ of the law of $X_t$ defined by \eqref{eqannealedbis}, which is given by a time change of the convolution product \eqref{eqconvoldens}. Let $\hat \mu_t$ the density at time $t$ of the  law of $Y_t$ defined by \eqref{eqannealed}, i.e. the solution of the annealed Langevin dynamics driven by the score $\nabla \ln \hat p_t$. Consider $\lambda_t$ such that $\lambda_0=0$.
    If $\nu$ is strictly log-concave and $\nabla\ln {p}_t$ is $L(t)$ Lipschitz with $L(t)\leq c$ for $t\in[0,T]$, then
        \begin{align}\label{eq:lsi_kl}
        d_{KL}(\tilde{\pi},\mathcal{L}(Y_{T/\kappa}))  & \leq \frac{\kappa}{2}  \int_0^{T}\vert p'_{s}\vert^2 e^{-\frac{1}{\kappa}\int_s^{T} C^{-1}_{\text{LS}}({\hat \mu}_{u/\kappa})d u} \, d s.
    \end{align}
\end{theorem}
Note that the integral term is expected to be small since $\kappa$ is small.

To understand what we intend to do, we will first give a formal derivation, assuming that there are no technical problems. We will then show how to overcome all the (unfortunately existing) technical problems.

\begin{proof}[Formal proof] 
Expanding the time derivative of the Kullback–Leibler distance, we have
    \begin{align*}
        \partial_t \;d_{KL}(\hat{p}_t, \hat{\mu}_t) & = \partial_t \int \ln \frac{\hat{p}_t}{\hat{\mu}_t} d \hat{p}_t = \int  \ln \frac{\hat{p}_t}{\hat{\mu}_t} \partial_t \hat{p}_t -\int\frac{\hat{p}_t}{\hat{\mu}_t} \partial_t \hat{\mu}_t \\
        & = -\int \ln \frac{\hat{p}_t}{\hat{\mu}_t} \, \nabla.(v_t\hat{p}_t) \, dx \, + \int\frac{\hat{p}_t}{\hat{\mu}_t} \nabla. \left(\hat{\mu}_t\nabla \left(\ln\frac{\hat{p}_t}{\hat{\mu}_t}\right)\right) \, dx\\
        &=\int \left\langle v_t, \nabla \ln \frac{\hat{p}_t}{\hat{\mu}_t}\right\rangle \hat{p}_t \, dx \, \, - \int \left\Vert \nabla \ln\frac{\hat{p}_t}{\hat{\mu}_t}\right\Vert^2 \hat{p}_t \, dx\\
        &\leq - \, (1-\frac{\lambda}{2})\int \left\Vert \nabla \ln\frac{\hat{p}_t}{\hat{\mu}_t}\right\Vert^2 \hat{p}_t \, dx + \frac{1}{2\lambda}\int \Vert v_t\Vert^2 \hat{p}_t \, dx\\
        &= - \, (1-\frac{\lambda}{2}) \int \left\Vert \nabla \ln\frac{\hat{p}_t}{\hat{\mu}_t}\right\Vert^2 \hat{p}_t  \, dx \, + \frac{\kappa^2}{2 \, \lambda} \vert p'_{\kappa t}\vert^2
    \end{align*}
for all $\lambda>0$. Recall that we have selected the Borel vector field $v_t$ such that it minimizes the $\mathbb{L}^2(\hat{p}_t)$-norm, i.e. $\Vert v_t\Vert_{\mathbb{L}^2(\hat{p}_td x)} = \kappa \vert p'_{\kappa t}\vert$. Notice that for $\lambda=2$ we recover \eqref{eqfinalKL}, Theorem \ref{eqannealed}. 
\medskip  
   
Assume that $\mu_t$ satisfies a logarithmic Sobolev inequality (this will be proved later). It follows that 
    \begin{align*}
        - \int \left\Vert \nabla \ln \frac{\hat{p}_t}{\hat{\mu}_t}\right\Vert^2 \hat{p}_t \, dx \, = \,  - \int \left\Vert \nabla \ln \frac{\hat{p}_t}{\hat{\mu}_t}\right\Vert^2 \frac{\hat{p}_t}{\hat \mu_t} \, \hat \mu_t \, dx\leq - \frac{2}{C_{\text{LS}}(\hat{\mu}_t)} d_{KL}(\hat{p}_t, \hat{\mu}_t).
    \end{align*}
    Substituting this with $\lambda=1$ we obtain
    \begin{align*}
        \partial_t \;d_{KL}(\hat{p}_t, \hat{\mu}_t) & \leq -\frac{1}{C_{\text{LS}}(\hat{\mu}_t)} d_{KL}(\hat{p}_t, \hat{\mu}_t) + \frac{\kappa^2}{2} \vert p'_{\kappa t}\vert^2.
    \end{align*}
    Applying Gronwall's inequality, we have
    \begin{align*}
        d_{KL}(\hat{p}_t, \hat{\mu}_t) & \leq \frac{\kappa^2}{2}  \int_0^{t}\vert p'_{\kappa s}\vert^2 e^{-\int_s^{t} C^{-1}_{\text{LS}}(\hat{\mu}_u)d u} \, d s.
    \end{align*}
    Therefore, for $t=T/\kappa$
    \begin{align*}
        d_{KL}(\tilde{\pi},\mathcal{L}(Y_{T/\kappa})) & \leq \frac{\kappa^2}{2}  \int_0^{T/\kappa}\vert p'_{\kappa s}\vert^2 e^{-\int_s^{T/\kappa} C^{-1}_{\text{LS}}(\hat{\mu}_u)d u} \, d s \\  & = \frac{\kappa}{2}  \int_0^{T}\vert p'_{s}\vert^2 e^{-\frac{1}{\kappa}\int_s^{T} C^{-1}_{\text{LS}}({\hat \mu}_{u/\kappa})d u} \, d s,
    \end{align*}
    where we have used that $\mathcal{L}(X_0) = \mathcal{L}(Y_0)$. 
\end{proof}
The proof of Theorem \ref{thmLSimprov} consists of two main components: establishing that $\mu_t$ satisfies a log-Sobolev inequality and justifying the formal computation of the time derivative of the relative entropy. 

The following result regarding functional inequalities for drifted Brownian motion will play a key role in showing that  $\mu_t$ satisfies a log-Sobolev inequality. 
\begin{lemma}\label{lemcurved}
Let $Z_t=Z_0 + B_t + \int_0^t \, b(s,Z_s) \, ds$ where $b(t,.)$ is $L(t)$-Lipschitz for $t \in [0,T]$. Then $$C_{LS}(Z_t) \leq  e^{\int_0^t \, L(u) \, du} \, C_{LS}(Z_0) \, +  \; \int_0^t \, e^{\int_0^s \, L(u) \, du} \,ds \,  .$$
Alternatively, if $b(t, .)$ satisfies $\nabla b(t, .)\leq - K(t) Id$, with $K(t)>0$ for all $t\in[0, T]$. Then
$$C_{LS}(Z_t) \leq  e^{-\int_0^t \, K(u) \, du} \, C_{LS}(Z_0) \, +  \; \int_0^t \, e^{-\int_0^s \, K(u) \, du} \,ds \,  .$$
\end{lemma}
This result is exactly Theorem 5 in \cite{CGsemin} since $$\langle b(t,x)-b(t,y),x-y\rangle \geq -L (t) \, |x-y|^2 \quad \text{or}\quad \langle b(t,x)-b(t,y),x-y\rangle \leq - K (t) \, |x-y|^2$$ for all $t$, which is denoted (H.C.-L(t)) or (H.C.K(t)) in \cite{CGsemin}. The proof, given in detail for the homogeneous case in Proposition 3 of \cite{CGsemin}, does not require any additional regularity for $b$. The only difference is our normalization of the log-Sobolev constant ($2 C_{LS}$ here is $C_{LS}$ there).

Before providing the proof of Theorem \ref{thmLSimprov}, we present conditions under which the bound in Theorem \ref{thmLSimprov} improves upon \eqref{eqfinalKL}.
\begin{proposition}\label{propconditionsimprovLS1}
    Let $\kappa\in(0,1/2)$, $\alpha\in(0,1/2]$ such that $\kappa^\alpha<1/2$, $\nu$ be a $\mathcal{N}(0, \sigma^2 I)$ and $\pi$ be compactly supported in the euclidean ball $B(0, R)$. Choose 
    \begin{gather*}
     \lambda(t) = 2(t/T)^2\mathbf{1}_{t< T/2} + (1-\kappa^{\alpha}-2\left(1-2\kappa^{\alpha}\right)(1-t/T)^{2})\mathbf{1}_{t\geq T/2}.
    \end{gather*}
    and $\sigma^2 = \frac{R^2}{\kappa^\alpha}$.
    Then, there exists a constant $C:=C(R, T)$ such that
\begin{equation*}
 d_{KL}(\tilde\pi,\mathcal{L}(Y_{T/\kappa})) \leq \, C\,d  \, \kappa^{2(1-\alpha)}.
 \end{equation*}
\end{proposition}
\begin{proof}
With our choice 
$$|p'_t|^2 \leq \frac{|\lambda'_t|^2}{2} \, \left(\frac{V_\pi}{\lambda_t} + \frac{V_\nu}{1-\lambda_t}\right) \leq \,\frac{8}{T^2} \, (R^2 + \sigma^2\, d)\,\leq \, 16\,\frac{R^2}{T^2\kappa^{\alpha}}\,d \,.$$
According to Proposition \ref{propgausscomp}, we have that for all $t\in[0, T]$
    $$
    \nabla^2\ln p_t \, \leq \, - \, \frac{\kappa^\alpha(1-\lambda_T)-\kappa^{2\alpha}\lambda_T}{R^2 \, (1-\lambda_T)^2} \, Id \, = -\frac{\kappa^\alpha}{R^2}\, Id\,.
    $$
Therefore, $p_t$ is strictly log-concave and applying Lemma \ref{lemcurved} (Theorem 5 in \cite{CGsemin}), it follows that
\begin{align*}
C_{LS}(Y_t) &\leq C_{LS}(\nu) e^{-T\kappa^\alpha/R^2} + 
\frac{4R^2(1-e^{-T\kappa^\alpha/R^2})}{\kappa^\alpha}
\leq\sigma^2 + \frac{4 R^2}{\kappa^\alpha} = \frac{5 R^2}{\kappa^\alpha}\,.
\end{align*}
Substituting all this into \eqref{eq:lsi_kl}, we obtain
\begin{equation*}
 d_{KL}(\tilde\pi,\mathcal{L}(Y_{T/\kappa})) \leq \, C(R, T)\,d \, \kappa^{2(1-\alpha)}\,.
 \end{equation*}
\end{proof}

Furthermore, if $\pi$ can be expressed as the convolution of a gaussian with a compactly supported distribution, an alternative result can be established.
\begin{proposition}\label{propconditionsimprovLS2}
    Let $\kappa\in(0,1/2)$, $\nu$ be a $\mathcal{N}(0, \sigma^2 I)$ and $\pi$ be the convolution of $\mathcal{N}(0, \tau^2 I)$ with a distribution supported in the euclidean ball $B(0, R)$, with $\tau^2\geq R^2$. Choose 
    \begin{gather*}
     \lambda(t) = 2\left(t/T\right)^2\mathbf{1}_{t\leq T/2} + \left(1-2(1-t/T)^2\right)\mathbf{1}_{t> T/2}.
    \end{gather*}
    Then, there exists a constant $C:=C(\tau^2, \sigma^2, R, T)$ such that
\begin{equation*}
 d_{KL}(\pi,\mathcal{L}(Y_{T/\kappa})) \leq \,C\, d\,\kappa^2.
 \end{equation*}
\end{proposition}
\begin{proof}
With our choice 
$$|p'_t|^2 \leq \frac{|\lambda'_t|^2}{2} \, \left(\frac{V_\pi}{\lambda_t} + \frac{V_\nu}{1-\lambda_t}\right) \leq \frac{4}{T^2}\,\left(R^2 + (\tau^2 + \sigma^2)\,d\right) .$$
On the other hand, according to Proposition \ref{propgausscompgauss}, we have that for all $t\in[0, T]$
    $$
    \nabla^2\ln p_t \, \leq \, - \, \frac{\min(\sigma^2, \tau^2-R^2)}{\max(\sigma^2, \tau^2)} \, Id \, = -K \, Id \,,
    $$
    where $K\geq0$.
Therefore, $p_t$ is strictly log-concave and applying Lemma \ref{lemcurved} (Theorem 5 in \cite{CGsemin}), it follows that
\begin{align*}
C_{LS}(Y_t) &\leq C_{LS}(\nu) e^{-KT} + \frac{4(1-e^{-KT})}{K} \leq \min\left(\sigma^2 + 4T, \sigma^2 + \frac{4}{K}\right) \,.
\end{align*}
For $K=0$, the second term in the first expression is $T$ instead of $(1-e^{-KT})/K$, however, the final bound still holds.
Plugging this into \eqref{eq:lsi_kl}, we obtain
\begin{equation*}
 d_{KL}(\pi,\mathcal{L}(Y_{T/\kappa})) \leq \, \frac{2}{T^2}\,\left(R^2 + (\tau^2 + \sigma^2)\,d\right) \, \min\left(\sigma^2 + 4T, \sigma^2 + \frac{4}{K}\right)  \, \kappa^2.
 \end{equation*}
 Note that we can take $\tilde\pi = \pi$ since $\lambda_T =1$.
\end{proof}

\begin{remark}
As shown in the proof below, thanks to Lemma \ref{lemcurved}, $\mathcal{L}(Y_{T/\kappa})$ satisfies a logarithmic Sobolev inequality. 
Therefore, we can use the implication between the log Sobolev inequality and the Talagrand $T_2$ inequality, Theorem 22.17 in \cite{Voldnew} to derive a bound in Wasserstein-2 distance
$$
W_2^2(\tilde\pi, \mathcal{L}(Y_{T/\kappa}))\leq 2 C_{LS}(Y_{T/\kappa})\,d_{KL}(\tilde\pi,\mathcal{L}(Y_{T/\kappa}))\,.
$$
In particular, under the assumptions of Propositions \ref{propconditionsimprovLS1} and \ref{propconditionsimprovLS2}, we have
$$
W_2^2(\tilde\pi, \mathcal{L}(Y_{T/\kappa})) \leq C(R, T) d \kappa^{2(1-3/2\alpha)}
\quad \text{and}\quad
W_2^2(\pi, \mathcal{L}(Y_{T/\kappa})) \leq C(\tau^2, \sigma^2, R, T) d \kappa^{2},
$$
respectively.
\hfill $\diamondsuit$
\end{remark}

\begin{proof}[Proof of Theorem \ref{thmLSimprov}]
First, we note that, under the  assumptions of the theorem, $\mathcal L(\sqrt{\lambda_0} \, X + \sqrt{1-\lambda_0} \, Z) = \mathcal{L}(Z)$ satisfies a log-Sobolev inequality and $L(t)$ is uniformly bounded on $[0, T]$. Therefore, Lemma \ref{lemcurved} provides control over the log-Sobolev constant for $\mu_t$ (uniformly in $t$).

Next, it remains to justify the formal calculation we have made for the time derivative of the relative entropy. To this end we will first collect some properties of $\hat p_t$ and $\hat \mu_t$. 

First $\hat p_t$ is the law of $X_t=\sqrt{\lambda_t} X + \sqrt{1-\lambda_t} Z$ so that $$V_{\hat p_t} \leq 2 (\lambda_t V_\pi + (1-\lambda_t) V_\nu) \leq 2 (V_\pi+V_\nu) \, .$$ Next $t \mapsto \nabla \ln(\hat p_t)$ is continuous, simply using Lebesgue's continuity theorem for integrals. It follows that $\nabla \ln(\hat p_t)(0)$ is bounded w.r.t $t\in[0,T]$. Since $\nabla \ln(\hat p_t)$ is $c$ global Lipschitz, and thanks to what precedes, it is at most of linear growth, uniformly in $t$, so that $$\sup_t \, \int \,  |\nabla \ln(\hat p_t)|^2 \, \hat p_t \, dx \, < \, +\infty \, .$$ Also notice that, since $e^{-W}$ (or $e^{-U}$) is bounded, $\hat p_t \in \mathbb L^1(dx) \cup \mathbb L^\infty(dx)$ hence belongs to all the $\mathbb L^q$ by interpolation.

We turn to the properties of $\hat \mu_t$. Using that $Y_.$ is a drifted Brownian motion with a $c$-Lipschitz drift, $Y_t$ has bounded (w.r.t. $t\in [0,T]$) moments of order $k$ as soon as $Y_0$ has a finite moment of order $k$. Since $\nu$ is strictly log-concave, $Y_0$ has finite moments of any order. In particular $$ \sup_t \, \mathbb E^{Q_Y}[|\nabla \ln \hat p_t(\omega_t)|^2] = \sup_t \, \int \,  |\nabla \ln \hat p_t|^2 \, \hat \mu_t \, dx \, < \, + \infty \, .$$ Let $B^\nu$ be the law of a Brownian motion with initial distribution $\nu$ (and variance $2$), denoted by $B_t^\nu$. The previous inequality shows that $$d_{KL}(Q_Y,B^\nu) \, < \, +\infty \, .$$ Denote by $R$ the time reversal operator at time $T$. Since the Kullback-Leibler distance decays using measurable push-forward, $d_{KL}(Q_Y\circ R,B_\nu\circ R) \leq d_{KL}(Q_Y,B_\nu) < + \infty$. 

Let $\gamma_t$ be the density of the law $\mathcal L(Z+\sqrt {2t} G)$ where $G$ is a standard gaussian measure. $\gamma_t$ is thus the density of $B_t^\nu$. All what is done in section \ref{secstoc} and section \ref{secgame} is unchanged if we replace $X_0$ by $G$, $\sqrt{\lambda_t}$ by $\sqrt{2t}$ and $\sqrt{1-\lambda_t}$ by $1$. We may thus again apply Theorem \ref{thmpoincmut} and  Lemma \ref{lemmalower} to deduce that $\nabla \ln(\gamma_t)$ is $c'$ global Lipschitz uniformly in $t\in [0,T]$. In particular, as for $\nabla \ln(\hat p_t)$, it is at most of linear growth uniformly in $t$. 

Standard results on time reversal thus say that $B^\nu\circ R$ is the law of a Brownian motion with drift $- \, 2 \nabla \ln(\gamma_{T-t})$ and initial measure $\gamma_T dx$. $Q_Y\circ R$ is thus the law of a Brownian motion with drift $r_{T-t}= - 2 \, \nabla \ln(\gamma_{T-t}) + h_{T-t}$ according to the Girsanov transform theory, and $$\int_0^T \, \int \mathbb E^{Q_Y}[|h_{T-t}\circ R|^2] \, dt < + \infty$$ thanks to the finite entropy condition (implying in particular $d_{KL}(\mathcal L(Y_T),\mathcal L(B_T^\nu))<+\infty$). 

Using the results in \cite{CCGL} (see Theorem 4.9 and formula (4.12) therein, also see the unpublished preliminary version \cite{CatPet} Corollary 3.15 and Theorem 4.8) one deduces the so called duality equation $$\nabla \ln (\hat \mu_t)=  \,  \nabla \ln(\hat p_{t}) \, + \, \bar h_{t}$$ for $t>0$, where $\bar h_{t}(x) = \mathbb E^{Q_Y}[h_{T-t}\circ R|X_t=x]$. In particular, using Cauchy-Schwarz inequality for the conditional expectation, $$\int_0^T \, \int |\bar h_{t}|^2 \, \hat \mu_t(dx) \,  \, dt \, < \, + \infty \, .$$ Since $\nabla \ln(\hat p_t)$ is of linear growth, it is square integrable w.r.t. $\hat \mu_t$ uniformly in $t$. All this implies that $$\int_0^T \, \int \, |\nabla \ln (\hat \mu_t)|^2 \, \hat \mu_t \, dx \, dt \, < \, + \infty \, .$$ Unfortunately this is not enough to ensure that $\nabla \ln(\hat p_t/\hat \mu_t) \in \mathbb L^1([0,T],\mathbb L^2(\hat p_t dx))$.
\medskip

We will thus have to use a (space) smooth cutoff function $\psi$ satisfying $\mathbf 1_{B(0,r)} \leq \psi \leq \mathbf 1_{B(0,R)}$ for some well chosen $0<r<R<+\infty$. Recall that $Y_.$ is a drifted Brownian motion with a drift $b(t,x)$ which is $c$-Lipschitz uniformly in $t$ and at most of linear growth, uniformly in $t$. One can thus use (a simplified and less precise form of) Theorem 3.1 in \cite{QZ} saying that the transition density of the process satisfies for all $T\geq t \geq t_0>0$, 
\begin{eqnarray*}
\hat \mu_t(x,y) &\geq& C(t_0,T) \, \exp \, - \, \left(c_1(t_0,T) |y-x|^2 \, e^{2\alpha t} \, + c_2(t_0,T) \sum_j |x_j-y_j|\right) \\ &\geq& C(t_0,T) \, \exp \, - \, \left(2c_1(t_0,T) |x e^{\alpha t}|^2 + c_2(t_0,T) \sum_j |x_j|\right) \, \\ && \quad \quad \quad \quad \exp \, - \, \left(2c_1(t_0,T) |y|^2 + c_2(t_0,T) \sum_j |y_j|\right),
\end{eqnarray*}
 for some $\alpha >0$, where all the constants depend on $C_W,D_W,C_U$ and $\sup_t |\nabla \ln(\hat p_t)(0)|$.

Actually this Theorem is only stated for a time homogeneous drift. A careful reading of the proof shows that it extends to a time inhomogeneous drift, provided it is of at most linear growth, uniformly in $t$. Indeed the proof is based on a comparison result with two other, time homogeneous, diffusion processes (Theorem 2.7) first shown in \cite{QRZ}. This Theorem immediately extends to time inhomogeneous drifts of at most linear growth uniformly in $t$ as conditions 1 and 2 in the Theorem are preserved. It is worth highlighting that the preliminary Corollary 2.6 in \cite{QZ} is written for a time dependent perturbation, and furnishes the proof of Theorem 2.7. 

We immediately deduce a rough lower bound
\begin{equation}\label{eqrough}
\hat \mu_t(x) \, \geq \, C'(t_0,T) \, e^{-c'_1(t_0,T) |x|^2 - c'_2(t_0,T) |x|} \, .
\end{equation}

\begin{remark}\label{remrough}
Estimates for general heat kernels have a long history, starting with Nash and overall Aronson in 1967, for generators in divergence form with very weak regularity. For non divergence generators a lot of work has been done for regular coefficients. For only Lipschitz, but unbounded, drifts the literature is not as rich. In addition to \cite{QRZ,QZ}, we may mention the more recent \cite{TaTa} covering more general situations (including non constant diffusion coefficient). We also refer to the Bibliography of the latter. 
\hfill $\diamondsuit$
\end{remark}
\medskip

Now consider $t\geq t_0>0$. Recall that $\psi$ is smooth and satisfies $\mathbf 1_{B(0,r)} \leq \psi \leq \mathbf 1_{B(0,R)}$ for some well chosen $0<r<R<+\infty$.

First
\begin{eqnarray*}
\int \, \psi \, \ln\frac{\hat{p}_t}{\hat{\mu}_t} \hat p_t \, dx - \int \, \psi \, \ln\frac{\hat{p}_s}{\hat{\mu}_s} \hat p_s \, dx &=& \int_s^t \, \int \, \psi \, (1+\ln \frac{\hat p_u}{\hat \mu_u}) \partial_u \hat p_u \, du - \, \int_s^t \, \int \, \psi \, \frac{\hat p_u}{\hat \mu_u} \, \partial_u \, \hat \mu_u \, du\\ &=& - \int_s^t \, \int \, \psi \, (1+\ln \frac{\hat p_u}{\hat \mu_u}) \nabla.(v_u \hat p_u) \, dx \, du \\ && \quad + \, \int_s^t \, \int \, \psi \, \frac{\hat p_u}{\hat \mu_u} \, \nabla.\left(\hat \mu_u \nabla \left( \ln \frac{\hat p_u}{\hat \mu_u}\right)\right) \, dx \, du \\ &=& \int_s^t \, \int \, \psi \, \left\langle v_u, \nabla \ln \frac{\hat{p}_u}{\hat{\mu}_u}\right\rangle \hat{p}_u \, dx \, du \\ && \quad - \, \int_s^t \, \int \, \psi \, ||\nabla \ln \frac{\hat p_u}{\hat \mu_u}||^2 \, \hat p_u \, dx \, du \\ && \quad + \, \int_s^t \, \int \, \langle \nabla \psi,v_u\rangle \, (1+\ln \frac{\hat p_u}{\hat \mu_u}) \, \hat p_u \, dx \, du \\ && \quad - \, \int_s^t \, \int \, \langle \nabla \psi,\nabla \ln \frac{\hat p_u}{\hat \mu_u}\rangle \, \hat p_u \, dx \, du \\ &=& A_1-A_2+A_3-A_4 \, .
\end{eqnarray*}
Everything is justified since
\begin{enumerate}
\item[(1)] \quad $v_u$ and $\nabla \ln \hat p_u$ belong to $\mathbb L^2(\hat p_u \, dx)$,
\item[(2)] \quad $\sqrt \psi \, \nabla \ln \hat \mu_u$ belongs to $\mathbb L^2(\hat p_u \, dx)$, since $\hat \mu_u$ is bounded from below by a positive constant on the support of $\psi$, hence $\psi \, \nabla \ln \frac{\hat p_u}{\hat \mu_u}$ belongs to $\mathbb L^2(\hat p_u \, dx)$ 
\item[(3)] \quad $\ln \hat p_u$ belongs to $\mathbb L^2(\hat p_u \, dx)$ since $\hat p_u$ belongs to all the $\mathbb L^q(dx)$, $\nabla \psi \, \ln \hat \mu_u$ is bounded, hence $\nabla \psi \, (1+\ln \frac{\hat p_u}{\hat \mu_u})$ belongs to $\mathbb L^2(\hat p_u \, dx)$,
\item[(4)] \quad $\langle\nabla \psi , \nabla \ln \frac{\hat p_u}{\hat \mu_u}\rangle$ belongs to $\mathbb L^2(\hat p_u \, dx)$ hence to $\mathbb L^1(\hat p_u \, dx)$, 
\end{enumerate}
and all the $\mathbb L^2$ norms in the previous items are integrable in time.

We now study each term $A_j$. First for $\lambda>0$, $$A_1 \leq \frac{\lambda}{2} \, \int_s^t \int \psi \left\Vert \nabla \ln\frac{\hat{p}_u}{\hat{\mu}_u}\right\Vert^2 \hat{p}_u \, dx \, du \, + \frac{1}{2\lambda} \int_s^t \, \int \Vert v_u\Vert^2 \hat{p}_u \, dx \, du \, .$$ Next we have $$A_3 \leq \frac{||\nabla \psi||_\infty}{2} \, \left(\int_s^t \, \int \mathbf 1_{|x|\geq r} \, ||v_u||^2 \, \hat p_u \, dx \, du \, + \, \int_s^t \, \int \mathbf 1_{|x|\geq r} \, (1+\ln \frac{\hat p_u}{\hat \mu_u})^2 \, \hat p_u \, dx \, du\right) \, ,$$ and 
\begin{eqnarray*}
\frac 13 \, \int_s^t \, \int \mathbf 1_{|x|\geq r} \, (1+\ln \frac{\hat p_u}{\hat \mu_u})^2 \, \hat p_u \, dx \, du &\leq&  \, \int_s^t \, \int \mathbf 1_{|x|\geq r} \, \hat p_u \, dx \, du \,  \, \\ && \quad \quad  +  \, \int_s^t \, \int \mathbf 1_{|x|\geq r} \, (\ln^2 \hat p_u) \, \hat p_u \, dx \, du   \\ && \quad \quad  \, + \, \int_s^t \, \int \mathbf 1_{|x|\geq r} \, (\ln^2 (1/\hat \mu_u)) \, \hat p_u \, dx \, du \, ,
\end{eqnarray*}
and finally thanks to \eqref{eqrough}, $$\int_s^t \, \int \mathbf 1_{|x|\geq r} \, (\ln^2 (1/\hat \mu_u)) \, \hat p_u \, dx \, du \, \leq \, C \, \int_s^t \, \int \mathbf 1_{|x|\geq r} \, (1+|x|^2)^2 \, \hat p_u \, dx \, du \, .$$ For $A_4$ we perform one more integration by parts. First we introduce a new cutoff $\eta$ similar to $\psi$ replacing $(r,R)$ by $(r/2,2R)$. This is done for the boundary term in the integration by parts to vanish. 

We thus have
\begin{eqnarray*}
A_4 &=&  \int_s^t \, \int \, \eta \, \langle \nabla \psi , \nabla \ln \frac{\hat p_u}{\hat \mu_u}\rangle \, \hat p_u \, dx \, du \, = \\ &=& - \, \int_s^t \, \int \, (\eta \, \hat p_u \, \Delta \psi + \langle \nabla \psi,\nabla(\eta \, \hat p_u)\rangle) \, \ln \frac{\hat p_u}{\hat \mu_u} \, dx \, du \\ &\leq& ||\Delta \psi||_\infty \, \int_s^t \, \int \mathbf 1_{|x|\geq r} \, \left|\ln \frac{\hat p_u}{\hat \mu_u}\right| \, \hat p_u \, dx \, du \, + \\ && \; +  \, ||\nabla \psi||_\infty \, \int_s^t \, \int \mathbf 1_{|x|\geq r} \, \left|\ln \frac{\hat p_u}{\hat \mu_u}\right| \, (||\nabla \eta||_\infty \, \hat p_u \, + \, \eta \, ||\nabla \hat p_u||) \, dx \, du \, = \, D_1 +D_2 \, .
\end{eqnarray*}
Gathering all this we have obtained, after simple manipulations,  $$(1-\frac{\lambda}{2}) \, \int_s^t \int \, \psi \left\Vert \nabla \ln\frac{\hat{p}_u}{\hat{\mu}_u}\right\Vert^2 \hat{p}_u \, dx \, du \\ \leq $$
\begin{eqnarray*}
&\leq& - \int \, \psi \, \ln\frac{\hat{p}_t}{\hat{\mu}_t} \hat p_t \, dx + \int \, \psi \, \ln\frac{\hat{p}_s}{\hat{\mu}_s} \hat p_s \, dx \, + \, \int_s^t \, \int \, \left(\frac{1}{2\lambda}+\frac{||\nabla \psi||_\infty }{2} \, \mathbf 1_{|x|\geq r}\right) \, \Vert v_u\Vert^2 \hat{p}_u \, dx \, du \, \\ && + \, \frac{3 \, ||\nabla \psi||_\infty}{2} \, \int_s^t \, \int \mathbf 1_{|x|\geq r} \, \left( 1 + \ln^2 \hat p_u + C (1+|x|^2)^2\right) \, \hat p_u \, dx \, du \, \\ && + \, \int_s^t \, \int \mathbf 1_{|x|\geq r} \, \left(||\Delta \psi||_\infty + ||\nabla \psi||_\infty ||\nabla \eta||_\infty\right) \, \left|\ln \frac{\hat p_u}{\hat \mu_u}\right| \, \hat p_u \, dx \, du \\ && + \, ||\nabla \psi||_\infty \, \int_s^t \, \int \mathbf 1_{|x|\geq r} \, \left|\ln \frac{\hat p_u}{\hat \mu_u}\right| \, ||\nabla \hat p_u|| \, dx \, du \, .
\end{eqnarray*}
The goal now is to let $r$ go to infinity, i.e. $\psi$ go to $1$. Of course this can be done assuming that $\nabla \psi$, $\Delta \psi$ and $\nabla \eta$ are uniformly bounded independently of $r$. Look at the right hand side. The first three terms in the sum converge to some limit, the fourth and the fifth one's go to $0$ thanks to the integrability properties we recalled (and the fact that $\int \rho |\ln \rho|  dx \leq c + \int \rho \ln \rho \, dx$). 

For the last term we remark that 
$$
 \left|\ln \frac{\hat p_u}{\hat \mu_u}\right| \, ||\nabla \hat p_u|| \leq  \, |\ln \hat p_u| \, \frac{||\nabla \hat p_u||}{\hat p_u} \, \hat p_u \,  +  \, C \, (1+|x|^2) \, \frac{||\nabla \hat p_u||}{\hat p_u} \, \hat p_u \, .$$ Recalling that $\ln \hat p_u$, $\frac{\nabla \hat p_u}{\hat p_u}$ and  $|x|^2$ all belong to $\mathbb L^2(\hat p_u)$, uniformly in time, we deduce thanks to Cauchy Schwarz inequality that the previous function is in $\mathbb L^1(\hat p_u)$, so that the last term goes to $0$ as $r \to \infty$, uniformly in time.
\medskip

We have finally obtained that $\nabla \ln \frac{\hat p_u}{\hat \mu_u} \in \mathbb L^2(\hat p_u)$, and in the limit $r \to \infty$ a justification for the final (formal) bound for the time derivative of the Kullback-Leibler distance. 

The proof is however not complete since all these derivations are made for $t\geq t_0>0$. It only means that we have first to integrate between $t_0$ and $T$, pass to the limit and then check that the remaining bounds do not depend on $t_0$ and conclude. Actually  this final step is not necessary, since it only modifies the constants.
\end{proof}

\begin{appendix}
\section*{}\label{secannex}
We recall the result from \cite{CL1}, which is used to establish the existence and uniqueness of weak solutions for \eqref{eqannealedbis} in Section \ref{secstoc}.
\begin{theorem}\label{thmCL}
Assume that $P$ is a unique weak solution of the martingale problem associated to the generator $L_P(t)=\frac 12 a(t,.) \partial^2 + b(t,.) .\nabla$, where $a=\sigma.\sigma^*$ with initial measure $\theta_0$. Let $\rho$ be an admissible flow in the following sense
\begin{enumerate}
\item[(i)] \quad $\rho_.$ satisfies the $(B,C_b^\infty)$ weak forward equation i.e. for all $f\in C_b^{1,\infty}(\mathbb R^+\times \mathbb R^d)$ and all $0\leq s <t \leq T$, $$\int f(t,x)d\rho_t - \int f(s,x)d\rho_s = \int_s^t \int (\partial_u f+ L_P(u) f+ \langle B(u,x),\nabla_x f\rangle) \, d\rho_u \, du \, .$$
\item[(ii)] \quad $B$ is of finite $\rho$-energy, i.e. for all $T\geq t>0$, $$\int_0^t \int \langle B,aB\rangle(s,x) \rho_s(dx) ds \, < \, +\infty \, .$$
\item[(iii)] \quad $H(\rho_0|\theta_0):=\int \ln(\frac{d\rho_0}{d\theta_0}) d\rho_0 <+\infty$.
\end{enumerate}
Assume in addition that one of the following conditions is satisfied
\begin{enumerate}
\item[(1)] \quad $\sigma$ and $b$ are $C^{1,2,\alpha}$ i.e. twice differentiable in $x$, one time differentiable in $t$ with derivatives that are $\alpha$ locally H\"{o}lder continuous,
\item[(2)] \quad $\sigma$ and $b$ are locally H\"{o}lder continuous, and $a$ is uniformly elliptic i.e. $\inf_{0\leq t \leq T} \, a(t,.) \geq \kappa \, Id$ in the sense of quadratic forms for some $\kappa >0$,
\item[(3)] \quad If $q_t$ denotes the marginal distribution $P\circ X_t^{-1}$ of $P$ at time $t$, $\rho_t $ is absolutely continuous w.r.t. $q_t$ and $d\rho_t/dq_t$ is locally bounded, i.e. bounded on $[0,T]\times K$ for any compact set $K \subset \mathbb R^d$.
\end{enumerate}
Then $Q$ defined by  
\begin{equation*}
\frac{dQ}{dP}|\mathcal F_T= \frac{d\rho_0}{d\theta_0} \, \exp \left(\int_0^T \, \langle B(s,\omega_s) \, , \,dM_s\rangle  \, -  \, \frac 12 \, \int_0^T \, \langle B,aB\rangle(s,\omega_s) \, ds\right)
\end{equation*}
where $M_.$ is the (local) martingale part of $X_.$ under $P$, is a probability measure satisfying $$H(Q|P) = H(\rho_0|\theta_0) + \frac 12 \, \int_0^T \int \langle B,aB\rangle(s,x) \rho_s(dx) ds \, < \, +\infty \, ,$$  whose time marginals at time $t$ are exactly $\rho_t$. 

In addition, $Q$ is the unique weak solution (up to time $T$) of the stochastic differential equation associated to $L_P(s)+aB(s,x).\nabla$ with initial distribution $\rho_0$.
\end{theorem} 
\bigskip
\end{appendix}
\begin{acks}[Acknowledgments]
A.~Guillin is supported by the ANR-23-CE-40003, Conviviality and has benefited from a government grant managed by the Agence Nationale de la Recherche under the France 2030
investment plan ANR-23-EXMA-0001.
P.~Cordero-Encinar is supported by EPSRC through the Modern Statistics and Statistical Machine Learning (StatML) CDT programme, grant no. EP/S023151/1.
\end{acks}

\bibliographystyle{imsart-number} 
\bibliography{references}       

\end{document}